\documentclass[a4paper, reqno]{amsart}
\usepackage{amssymb, amsmath, amsthm, chngcntr, enumerate, mathabx, mathrsfs, mathtools,esint}
\usepackage{dsfont}
\usepackage[english]{babel}
\usepackage[T1]{fontenc}

\def\be#1{\begin{equation}\label{#1}}
\def\bas{\begin{align*}}
\def\eas{\end{align*}}
\def\bi{\begin{itemize}}
\def\ei{\end{itemize}}

\theoremstyle{plain}
   \newtheorem{theorem}[subsection]{Theorem}
   \newtheorem*{theorem*}{Theorem}

   \newtheorem{proposition}[subsection]{Proposition}

\theoremstyle{remark}
   \newtheorem{remark}{Remark}

\theoremstyle{definition}

\author[M. Vitturi]{Marco Vitturi}
\address{Marco Vitturi: CNRS - Universit\'{e} de Nantes \\ Laboratoire Jean
Leray \\ 2, rue de la Houssini\`{e}re
44322 Nantes cedex 3, France}
\email{marco.vitturi@univ-nantes.fr}

\author[J. Wright]{James Wright}
\address{James Wright: School of Mathematics and Maxwell Institute for Mathematical Sciences, University of
Edinburgh, JCMB, King's Buildings, Peter Guthrie Tait Road, Edinburgh EH9 3FD, Scotland}
\email{J.R.Wright@ed.ac.uk}


\title[Multiparameter singular integrals on the Heisenberg group]{Multiparameter singular integrals on the Heisenberg group: uniform estimates}

\begin{document}

\begin{abstract} We consider a class of multiparameter singular Radon integral operators
on the Heisenberg group ${\mathbb H}^1$ where the underlying variety is the graph of a polynomial. A remarkable difference
with the euclidean case, where Heisenberg convolution is replaced by euclidean convolution, is that
the operators on the Heisenberg group are always $L^2$ bounded. This is not the case
in the euclidean setting where $L^2$ boundedness depends on the polynomial defining
the underlying surface. Here we uncover some new, interesting phenomena. For example,
although the Heisenberg
group operators are always $L^2$ bounded, the bounds are {\it not} uniform in the coefficients of 
polynomials with fixed degree. When we ask for which polynoimals uniform $L^2$ bounds hold,
we arrive at the {\it same} class where uniform bounds hold in the euclidean case.
\end{abstract}

\maketitle

\section{Introduction}
For the general theory of singular Radon transforms
$$
H_{\gamma, K} f(x) \ = \ \psi(x) \int_{{\mathbb R}^k} f(\gamma(x,t)) K(t) \, dt
$$
where $K$ is a singular kernel and $\gamma : {\mathbb R}^n \times {\mathbb R}^k \to {\mathbb R}^n$
is a smooth map ($\psi$ an appropriate cut-off function), the case of
translation-invariant polynomial
mappings $\gamma(x,t) = x \cdot {\Phi}(t)$ has served as a model problem.
Here ${\Phi}(t) = (P_1(t), \ldots, P_n(t))$ with 
polynomial components $P_j \in {\mathbb R}[X_1,\ldots, X_k]$ and the translation $\cdot$ arises
from a nilpotent Lie group structure on ${\mathbb R}^n$. See \cite{CNSW} where the analysis of general
singular Radon transforms $H_{\gamma, K}$ is effectively reduced to the case 
$\gamma(x,t) = x \cdot {\Phi}(t)$
described above in the one-parameter setting; that is, when $K$ is a classical Calder\'on-Zygmund CZ
kernel; $|\partial^{\alpha} K (t)| \lesssim |t|^{-k-|\alpha|}$ with appropriate cancellation conditions
imposed.

In polynomial translation-invariant cases $\gamma(x,t) = x \cdot {\Phi}(t)$, one
heuristic idea crucial to the study of
the operator $H_{\gamma, K} = H_{{\Phi}}$ is that if the boundedness of 
$H_{{\Phi}}$ is to be proved, one
must do so by proving the stronger statement that the bound can be taken
to be independent of the polynomial ${\Phi}$, once the degree of ${\Phi}$ is fixed. 
This is especially the case in the one-parameter setting; 
see \cite{RS1} where this heuristic is developed systematically and consequences are 
explored.


For multiparameter CZ kernels $K$ (see Section \ref{street} for a precise definition), 
the operators $H_{\gamma, K}$ may or may not be
$L^2$ bounded and matters depend on cancellation conditions which
arise through a subtle interaction between the mapping $\gamma$ and the kernel $K$. In the euclidean translation-invariant setting, these cancellations conditions have been thoroughly investigated by Ricci and Stein in \cite{RS3} (see \cite{NW} for earlier work). In particular Theorem 5.1 in \cite{RS3} gives a sufficient condition (a cancellation condition involving both $\gamma$ and $K$) which guarantees $L^2$ (even $L^p$) boundedness of the associated singular integral operators. One can then check in particular instances if these conditions are necessary.  

For instance if $\gamma(x,t) = x + \Sigma(t)$ where
$\Sigma(t) = (t, P(t)) $ parametrises
an $(n-1)$-dimensional polynomial surface with $P \in {\mathbb R}[X_1, \ldots, X_{n-1}]$,
then the so-called Hilbert transform along $\Sigma$, 
$H_{\gamma, {\mathcal K}} = {\mathcal H}_{P, {\mathcal K}}$ where
\[
{\mathcal H}_{P, {\mathcal K}} f(x, z) \ = \ p.v. \int_{{\mathbb R}^{n-1}} f(x - t, z - P(t)) \,
{\mathcal K}(t) \, dt,
\]
is a typical example of a multiparameter singular Radon transform treated in
\cite{RS3} (see also \cite{NW}). Here the multiple Hilbert transform kernel
${\mathcal K}(t) = 1/t_1\cdots t_{n-1}$ is the canonical multiparameter CZ kernel.
If  $P(t) = \sum_{\alpha} c_{\alpha} t^{\alpha}$ is a real polynomial in $n-1$ variables,
we define the \emph{support of $P$} as $\Delta_P = \{\alpha: c_{\alpha} \neq 0\}$.
For any finite $\Delta \subset {\mathbb N}_0^{n-1}$, let 
${\mathcal V}_{\Delta}$ denote 
the finite dimensional subspace of real polynomials $P$ in $n$ variables with $\Delta_P \subseteq \Delta$.

\begin{theorem*}[Ricci-Stein \cite{RS3}]
Fix $\Delta \subset {\mathbb N}_0^{n-1}$. Then
\begin{equation}\label{RS-L2}
\sup_{P \in {\mathcal V}_{\Delta}} \|{\mathcal H}_{P, {\mathcal K}} \|_{L^2 \to L^2} \ < \ \infty
\end{equation}
holds if and only if every $\alpha = (\alpha_1, \ldots, \alpha_{n-1}) \in \Delta$, at least $n-2$
of the $\alpha_j$'s are even. Furthermore if $\alpha$ has 2 odd components, then for
$P(t) = t^{\alpha}$, the individual operator
${\mathcal H}_{P, {\mathcal K}}$ is unbounded on $L^2$.
\end{theorem*}

More precisely, the sufficiency part of this theorem follows from 
Theorem 5.1 in \cite{RS3} via a standard lifting procedure (effectively freeing up the
monomials of $P$) 
to an operator on a higher dimensional space of the form ${\mathcal H}_{Q,{\mathcal K}}$ where 
$$
Q(t) \ = \ (Q_{\alpha}(t) )_{\alpha \in \Delta_P} \ \ \ {\rm and \ each} \ \ 
Q_{\alpha}(t) = t^{\alpha}.
$$
One then checks that $Q$ and ${\mathcal K}$ satisfy the cancellation condition of Theorem 5.1 in
\cite{RS3}. For the necessity it is a simple computation to check that if $P(t) = t^{\alpha}$
and $\alpha$ has 2 odd components, then ${\mathcal H}_{P,{\mathcal K}}$ is unbounded on $L^2$ (see \cite{Fef}).

This result depends very much on the multiparameter CZ kernel under consideration. If the
multiple Hilbert transform kernel ${\mathcal K}$ is replaced by a different multiparameter CZ kernel,
the cancellation condition in Theorem 5.1 changes. See \cite{B2} where a {\it projected}
version of ${\mathcal H}_{P, K}$ is considered for a fixed polynomial $P$ but the multiparameter
CZ kernels $K$ varies. A sharp result is established where uniformity in $K$ is sought for a fixed
polynomial $P$.

In a remarkable series of papers, the translation-invariant theory of Ricci and Stein 
was extended to the general
non-translation-invariant setting by Stein and Street; \cite{SS1,B1,B2,SS2,SS3}
and \cite{B}. In this work two conditions on $\gamma$ are introduced, one is a curvature
condition generalising the fundamental curvature condition in \cite{CNSW} and another is
an algebraic condition which can be viewed as a strong cancellation condition. When these
two conditions hold, $L^2$ bounds for $H_{\gamma, K}$ are deduced for any multiparameter CZ kernel $K$.
These two conditions depend only on $\gamma$ and so the cancellation condition is decoupled from
the particular singular kernel under consideration. Hence the results obtained are valid for all
multiparameter CZ kernels. In many cases,
when uniformity in $K$ is sought, the algebraic or cancellation condition can be shown to be necessary.
See \cite{B2} for details.

A fascinating example is given by $\gamma({\underline{x}},s,t) = {\underline{x}} \cdot \Sigma(s,t)$ 
where $\Sigma(s,t) = (s,t,P(s,t))$
parameterises the graph of a polynomial surface in ${\mathbb R}^3$ and 
$\cdot$ is the Heisenberg group ${\mathbb H}^1 \simeq {\mathbb R}^3$ multiplication;
$(x,y,z)\cdot (u,v,w) = (x+u, y+v, z+w + 1/2(xv - yu))$. Interestingly, both conditions alluded to above
are always satisfied (see Section \ref{street} for details) and hence in particular, 
$H_{\gamma, {\mathcal K}}$ is bounded on
$L^2$ for {\it any} real polynomial $P$. This is in
sharp contrast to the above Ricci-Stein theorem which shows that in the euclidean translation-invariant case  
$\gamma({\underline{x}}, s,t) = {\underline{x}} + \Sigma(s,t)$, $L^2$ boundedness depends
on the particular polynomial $P(s,t)$. This extends to any real-analytic $P$ and any multiparameter
CZ kernel $K$. See \cite{B2} and \cite{SS2}.

\begin{theorem}\label{street-analytic} For any real polynomial $P(s,t)$ (or more generally
any real-analytic $P$ near the origin $(0,0)$) and multiparameter CZ kernel $K$, consider
$$
H_{P, K} f(x,y,z) \ = \ \iint_{{\mathcal R}} f( (x,y,z) \cdot (s,t,P(s,t))^{-1}) \,  K(s,t) \, ds dt
$$
where ${\mathcal R} = {\mathcal R}_{a,b,c,d} = \{(s,t): a\le |s| \le b, c \le |t| \le d\}$ 
is any ``rectangle'' but when $P$ is 
real-analytic at the origin, we take
$b$ and $d$ to be sufficiently small.
Then
$H_{P, K}$ is bounded on $L^2({\mathbb H}^1)$.
\end{theorem}

The arguments developed in this paper will not only give an alternative proof of Theorem \ref{street-analytic}
but will shed light on the r\^{o}le Heisenberg translations play in multiparameter settings. See Section \ref{street}
for an extension of Theorem \ref{street-analytic}. 
Interestingly when we seek $L^2$ bounds, uniform with respect to the polynomial $P$ as in
the Ricci-Stein theorem, we come back to the euclidean conclusion.

\begin{theorem}\label{main-uniform} Fix $\Delta \subset {\mathbb N}_0^{2}$. Then
\begin{equation}\label{RS-L2-Heisenberg}
\sup_{P \in {\mathcal V}_{\Delta}} \|{H}_{P, {\mathcal K}} \|_{L^2({\mathbb H}^1) \to L^2({\mathbb H}^1)} \ < \ \infty
\end{equation}
holds if and only if every $\alpha = (\alpha_1, \alpha_{2}) \in \Delta$ has at least one
even component.

More generally, for $H_{P, K}$ where $K$ is a general multiparameter CZ kernel $K$,
the uniformity in \eqref{RS-L2-Heisenberg} is equivalent to the uniformity of a family of truncations of the
singular Radon transform
\[
R_{P,K} g(x,y) \ = \ \iint\limits_{{\mathbb R}^2} g(x - t, y - P(s,t)) \, K(s,t) \, ds dt.
\]
\end{theorem} 

In this theorem the integration defining the operator $H_{P, {\mathcal K}}$ is taken over a ``rectangle''
${\mathcal R} = \{a\le |s| \le b, c \le |t| \le d\}$ and the uniformity conclusion holds with respect to the parameters
$a,b,c$ and $d$ as well.

{\bf Notation} Uniform bounds for oscillatory integrals lie at the heart of this paper. Keeping track of constants
and how they depend on the various parameters will be important for us. For the most part, constants $C$
appearing in inequalities $A \le C B$ between positive quantities $A$ and $B$ will be {\it absolute} or
{\it uniform} in that they can be taken to be independent of the parameters of the underlying problem. 
We will use $A \lesssim B$ to denote $A \le C B$ and $A \sim B$ to denote $C^{-1} B \le A \le C B$.
If $A$ is a general real or complex quanitity,
we write $A = O(B)$  to denote $|A| \le C B$ and when we want to highlight a dependency on a parameter
$\theta$, we write $A = O_{\theta}(B)$ or $|A| \lesssim_\theta B$ to denote $|A| \le C_{\theta} B$.


\section{The work of Street \cite{B2} and further results}\label{street}

In \cite{B2}, Street develops the
$L^2$ theory for multiparameter singular Radon transforms 
\[
H_{\gamma, K} f(x) \ = \ \psi(x) \int_{{\mathbb R}^k} f(\gamma(x,t)) \, K(t) \, dt
\]
and introduces two key conditions on $\gamma$; a finite-type (curvature) condition and an
algebraic (cancellation) condition. Here
$\gamma: {\mathbb R}^n \times {\mathbb R}^k \to {\mathbb R}^n$ is a smooth map satisfying
$\gamma(x,0) \equiv x$, $\psi$ an appropriate cut-off function, and $K(t)$ is multiparameter Calder\'on-Zygmund
kernel which is usually supported near the origin $t=0$.

For our purposes it suffices to restrict our attention to the
2-parameter case ${\mathbb R}^k = {\mathbb R}^{k_1} \times {\mathbb R}^{k_2}$ and to {\it product kernels}
$K$ as introduced in \cite{NRS}, which underpins the theory of singular integrals with respect to flag kernels (however 
our analysis extends to treat the more general class of multiparameter CZ kernels considered in \cite{B2}).

The notion of product kernel depends on the classical notion of CZ kernels in one parameter; that is, a distribution
$K$ on ${\mathbb R}^k$ which coincides with a smooth function away from the origin such that
$|\partial^{\alpha} K(t)| \le C_{\alpha} |t|^{k - |\alpha|}$ for all $\alpha$ and  such that the quantities
$\int K(t) \phi(R t) dt$ are bounded, uniformly over all $R>0$ and all smooth $\phi$ supported in the unit ball
with $\|\phi\|_{C^1} \le 1$ (such a $\phi$ is called a {\it normalised bump function} on ${\mathbb R}^k$).

A 2-parameter product kernel $K$ is defined as follows. It is 
a distribution on ${\mathbb R}^k = {\mathbb R}^{k_1} \times {\mathbb R}^{k_2}$ which coincides with a
$C^{\infty}$ function $K$ away from the coordinate subspaces $s = 0$, $t = 0$ and satisfies

1. (Differential inequalities) for every multi-index $\alpha = (\alpha_1, \alpha_2) \in \mathbb{N}^{k_1} \times \mathbb{N}^{k_2}$, there is a constant $C_{\alpha}$ such that
\[
| \partial^{k_1}_{\alpha_1} \partial^{k_2}_{\alpha_2} K (s,t)| \ \le \ C_{\alpha} \, |s|^{k_1 - |\alpha_1|} 
|t|^{k_2 - |\alpha_2|}
\]
away from the two coordinate subspaces, and

2. (Cancellation) for any normalised bump function $\phi$ on ${\mathbb R}^{k_1}$ 
and any $R>0$, the distribution
\[
K^1_{\phi, R} (t) \ = \ \int_{{\mathbb R}^{k_1}} K(s, t) \phi(R s) \, ds
\]
is a classical one parameter CZ kernel on ${\mathbb R}^{k_2}$ described above. Similarly for
$K^2_{\phi, R} (s) = \int K(s,t) \phi(R t) dt$. 

Important for our analysis is the following characterisation of product kernels; see Corollary 2.2.2
in \cite{NRS}. For every smooth $\phi$ and $I = (j,k) \in {\mathbb Z}^2$, we set 
$\phi^{(I)}(s,t) := 2^{-j-k} \phi(2^{-j} s, 2^{-k} t)$.

\begin{proposition}\label{product-characterisation} A product kernel $K$ can be written as
\begin{equation}\label{product-sum}
K \ = \ \sum_{I\in {\mathbb Z}^2} \phi_I^{(I)}
\end{equation}
(which is convergent in the sense of distributions) where each smooth $\phi_I$ is supported in
$\{(s,t) : 1/2\le |s|, |t|\le 2\}$, satisfies the cancellation conditions
\begin{equation}\label{product-cancellation}
\int \phi_I(s,t) \, ds \ \equiv \ 0 \ \ {\rm and} \ \ \int \phi_I(s,t) \, dt \ \equiv \ 0
\end{equation}
for every $t$ and $s$, and the sequence $\{\phi_I\}$ is bounded in $C^k$ norm for every $k$.
\end{proposition}

The two key conditions on $\gamma$ are easily formulated in the case
where $\gamma$ can be written as the exponential\footnote{The multiparameter exponential is to be interpreted as follows: for $s,t$ given, define vector field $Y_{s,t}(x)=\sum s^p t^q X_{p,q}(x)$; then $\exp\big(\sum s^p t^q X_{p,q}\big)(x) := \exp(\tau Y_{s,t})\big|_{\tau=1}(x)$. }
\begin{equation}\label{gamma_exponential}
\gamma(x,(s,t)) \ = \ \exp\big(\sum s^p t^q X_{p,q}\big) (x)
\end{equation}
of a \emph{finite} sum of smooth vector fields $\{X_{p,q} = X_{\alpha}\}$. We assign to each $X_{\alpha}$, where $\alpha = (\alpha_1,\alpha_2) \in \mathbb{N}^{k_1} \times \mathbb{N}^{k_2}$, the
formal degree $d_\alpha = (|\alpha_1|,|\alpha_2|) \in \mathbb{N} \times \mathbb{N}$ and recursively we then define formal degrees for all iterated commutators
such that if $d_1$ and $d_2 \in {\mathbb N}^2$ are the degrees of iterated commutators
$X_1$ and $X_2$, respectively, then $[X_1, X_2]$ has degree $d_1 + d_2 \in \mathbb{N} \times \mathbb{N}$. Hence we
view these vector fields, together with their corresponding degree $(X,d)$. Notice that it might be the case that one vector field has more than one degree; in this case we consider them to be distinct objects.

We separate the original vector
fields $\{(X_{\alpha}, d_\alpha)\} = {\mathcal P} \cup {\mathcal N}$ into two types; the pure ones 
$(X_{\alpha}, d_\alpha) \in {\mathcal P}$
where $d_\alpha = (p,0)$ or $d_\alpha=(0,q)$ and non-pure ones $(X_{\alpha}, d_\alpha) \in {\mathcal N}$ where 
$d_\alpha = (p,q)$ and both $p$ and $q$ are nonzero. The two key conditions on $\gamma$ are the following:
there is a finite list $\{(X_1, d_1), \ldots, (X_N, d_N)\}$ of iterated commutators of pure vector fields, containing ${\mathcal P}$ itself and such that

1. (Finite-type condition)  for all $\delta \in [0,1]^2$, we can write\footnote{Here $\delta^{d} = \delta_1^{d_1} \cdot \delta_2^{d_2}$.}
\begin{equation}\label{finite-type-curvature}
[\delta^{d_j} X_j, \delta^{d_k} X_k] \ = \ \sum_{\ell=1}^N c_{j,k}^{\ell, \delta} \delta^{d_{\ell}} X_{\ell}
\end{equation}
where $c_{j,k}^{\ell, \delta} \in C^{\infty}$, uniformly in $\delta$; and

2. (Algebraic condition) for $(Y, e) \in {\mathcal N}$ and every $\delta \in [0,1]^2$, we can write
\begin{equation}\label{algebraic-cancellation}
\delta^{e} Y \ = \ \sum_{\ell=1}^N  c_{Y}^{\ell, \delta} \delta^{d_{\ell}} X_{\ell}
\end{equation}
where $c_{Y}^{\ell, \delta} \in C^{\infty}$, uniformly in $\delta$.

\begin{remark}
Notice the two conditions imply that the involutive distribution generated by the collection $\{X_\alpha\}$ is finitely generated (as a $C^\infty$-module). In the one-parameter case, this is essentially equivalent to the conditions above and the scaling factors in $\delta$ play essentially no active r\^{o}le. However, this is no longer necessarily true in the multiparameter case (see \cite{B2}, Section 17.7) and the uniform behaviour in $\delta$ becomes crucial.
\end{remark}

The finite-type condition \eqref{finite-type-curvature} is a generalisation of the curvature condition
introduced in \cite{CNSW} in the one-parameter setting and the algebraic condition \eqref{algebraic-cancellation}
allows us to {\it control} the troublesome non-pure vector fields $Y\in {\mathcal N}$ in terms of the
pure ones, effectively transferring any needed cancellation down to the product kernel $K$.
In this case, under these two conditions on $\gamma$, $L^2$ bounds for $H_{\gamma, K}$ can be derived
for any product kernel $K$. In more general (non-finite, that is when $\gamma$ is not exactly of type \eqref{gamma_exponential}) situations,
the conditions \eqref{finite-type-curvature} and \eqref{algebraic-cancellation} need to be modified.
See \cite{B2} for details and in particular see section 3 of \cite{B2} for a discussion of the finite case discussed above.

The particular situation we are concerned with here  is $\gamma({\underline{x}}, (s,t)) = {\underline{x}} \cdot \Sigma(s,t)$
where the product $\cdot$ is the Heisenberg ${\mathbb H}^1$ group multiplication and $\Sigma(s,t) = (P_1(s,t), P_2(s,t), P_3(s,t))$
parametrises a surface in ${\mathbb H}^1$. Let $X = \partial_x - (y/2) \partial_z,  Y = \partial_y + (x/2) \partial_z$
and $Z = \partial_z$ be the usual basis of left-invariant vector fields on ${\mathbb H}^1$ such that 
$[X,Y] = Z$. Then 
\[
\gamma({\underline{x}}, (s,t)) \ = \ {\underline{x}} \cdot \Sigma(s,t) \ = \ {\rm exp}({P_1(s,t) X + P_2(s,t) Y + P_3(s,t) Z})
({\underline{x}}),
\]
putting us in the above {\it finite} situation if each $P_j$ is a polynomial. In this case
the finite-type condition \eqref{finite-type-curvature} is automatically satisfied. In turns out
that when the $P_j$ are (more generally) real-analytic, the appropriately modified finite-type
condition \eqref{finite-type-curvature} is still automatically satisfied; see \cite{SS2}.

In the case that 
$P_1(s,t) = s$ and $P_2(s,t) = t$, we see that $(X, (1,0))$ and $(Y,(0,1))$ lie in ${\mathcal P}$. Furthermore the only vector fields lying
in ${\mathcal N}$ must be of the form $(Z,d)$ where $d = (p,q)$ satisfies $p q \neq 0$ and the monomial $s^p t^q$
arises in the Taylor expansion of $P_3(s,t)$. Hence for any real-analytic $P_3$,
every non-pure vector field in ${\mathcal N}$ can be
controlled as described in \eqref{algebraic-cancellation} and so both conditions
\eqref{finite-type-curvature} and \eqref{algebraic-cancellation} are automatically satisfied
when $\Sigma(s,t) = (s,t, P(s,t))$ is the graph of a real-analytic surface in ${\mathbb H}^1$. 
This is the background discussion for Theorem \ref{street-analytic}.

Now let us consider a slight variant; a surface parameterised by 
$\Sigma(s,t) = (s^{p_0}, t , P(s,t))$ where $P$ is a general real-analytic
function near $(0,0)$. As mentioned above, the corresponding finite-type condition \eqref{finite-type-curvature}
is automatically satisfied but now it is not necessarily the case that all non-pure vector fields $(Z,d') \in {\mathcal N}$
can be controlled by pure vector fields in the sense of \eqref{algebraic-cancellation}. 
Recall that $d' = (p',q')$ where $p' q' \neq 0$ and $s^{p'} t^{q'}$ arises in
the series expansion of $P(s,t) = \sum c_{p,q} s^p t^q$. Note that if $p'\ge p_0$, then we
can control $(Z,d')$ by $(Z, d_0)$ where $d_0 = (p_0,1)$ and $(Z, d_0)$ arises
as the commutator of the pure vector fields $(X, (p_0,0))$ and $(Y,(0,1))$. Therefore the
non-pure vector fields $(Z,d')$ which cannot be controlled in the sense of \eqref{algebraic-cancellation}
must necessarily satisfy $p' < p_0$ and so arise from a term in 
$$
P_{p_0}(s,t) \  := \ \sum_{p=0}^{p_0 -1} \sum_{q\ge 1} c_{p,q} s^p t^q \ + \ \sum_{p\ge 1} c_{p,0} s^p.
$$
When $p_0 =1$, we have $P_{p_0}(s,t) = P(s,0)$
and so no $d' = (p', q')$ with $p' q' \neq 0$ satisfies $p' < p_0 =1$,
bringing us back
to the case where all non-pure terms can be controlled by pure ones; that is, condition \eqref{algebraic-cancellation}
is satisfied.

\begin{theorem}\label{vw-analytic} For any real-analytic $P(s,t)$ near the origin $(0,0)$ and
multiparameter CZ kernel $K$, consider
\[
H_{P, K} f(x,y,z) \ = \ \iint\limits_{{\mathcal R}} f((x,y,z)\cdot(s^{p_0}, t, P(s,t))) \,  K(s,t) \, ds dt
\]
where ${\mathcal R} = {\mathcal R}_{a,b,c,d} = \{(s,t): a\le |s| \le b, c \le |t| \le d\}$ lies in a small
neighbourhood of the origin $(0,0)$. If $P_{p_0} \equiv 0$, then $H_{P,K}$ is bounded on $L^2({\mathbb H}^1)$.

In general, the $L^2({\mathbb H}^1)$ boundedness of $H_{P,K}$ is equivalent to the uniform 
$L^2({\mathbb R}^2)$ boundedness of
a family of truncations of the singular Radon transform
\[
R_{P_{p_0}, K} g(x,y) \ = \ \iint\limits_{{\mathbb R}^2} g(x - t, y - P_{p_0}(s,t)) K(s,t) \, ds dt .
\]

Furthermore when $K$
is the double Hilbert transform kernel ${\mathcal K}(s,t) = 1/st$, then 
$H_{P, {\mathcal K}}$ is bounded on $L^2({\mathbb H}^1)$ 
if and only if every vertex $(p,q)$ of the Newton polygon of
$P_{p_0}$ has the property that $p q$ is even.
\end{theorem}

The Newton polygon of $P_{p_0}$ is the convex hull of the quadrants $(p,q) + {\mathbb R}^2_{+}$
in ${\mathbb R}^2$ wihere $(p,q) \in \Delta(P_{p_0})$, the support of $P_{p_0}$. The
 r\^{o}le of the Newton polygon in the theory of multiparameter singular Radon transforms
 first appeared in \cite{CWW1}.

The first part of Theorem \ref{vw-analytic} follows from the work of Stein and Street, \cite{B2} and \cite{SS2}.
The more general statement gives a precise structural description of the $L^2$ boundedness properties for $H_{P,K}$
and highlights the r\^{o}le of Heisenberg translations in multiparameter settings. Theorem \ref{vw-analytic} is
a representative theorem and exposes a new phenomenon for multiparameter convolution operators on the
Heisenberg group. More general results can be formulated and established. 

\section{Initial reductions for the proofs of Theorems 
\ref{main-uniform} and \ref{vw-analytic}.}\label{prelim-reduct}
We fix a product kernel $K$ and 
use Proposition \ref{product-characterisation} to write $K = \sum_I \phi_I^{(I)}$ as in \eqref{product-sum}
with the smooth, compactly supported $\phi_I$ satisfying \eqref{product-cancellation}. We consider the operator
\[
T_{P, {\mathcal F}} f(x,y,z) \ = \ 
\sum_{I\in {\mathcal F}} \int f( (x,y,z) \cdot (s^{p_0}, t, P(s,t))) \phi^{(I)}_I (s,t) \, ds dt
\]
where ${\mathcal F} \subset {\mathbb Z}^2$ is a fixed finite subset ${\mathcal F} = \{ I = (j,k)\}$,
indexing the dyadic rectangles ${\mathcal R}_I = \{(s,t): |s| \sim 2^j, |t| \sim 2^k\}$ where $\phi_I^{(I)}$
(and hence the integral above) is supported. For Theorem \ref{vw-analytic}, when $P$ is 
assumed to be real-analytic near the origin, we require that the rectangles ${\mathcal R}_I$
are supported near the origin; that is, if $I = (j,k) \in {\mathcal F}$, then both $j$ and $k$ are sufficiently negative. 

By translation-invariance in the third variable we may assume, without loss of generality, that 
$P(0,0) = 0$. Furthermore, the structure of the Heisenberg group allows us to make another reduction that will be very useful in the following. If $P(s,t) = c s^{p_0} + d t + {\tilde P}(s,t)$ with 
$\partial^{p_0}_s {\tilde P}(0,0) = \partial_t {\tilde P}(0,0) = 0$,
then we can write
\[
f((x,y,z)\cdot (s^{p_0}, t, P(s,t))) = f(L [L^{-1}(x,y,z) \cdot (s^{p_0}, t, {\tilde P}(s,t))])
\]
where 
\[
L \ = \ 
\begin{pmatrix}
1 & 0 & 0 \\
0 & 1 & 0 \\
c & d & 1 
\end{pmatrix}
\] 
is a \emph{group automorphism} of ${\mathbb H}^1$.
Hence $\|T_{P, {\mathcal F}}\|_{L^2 \to L^2} = \|T_{{\tilde P}, {\mathcal F}}\|_{L^2 \to L^2}$ and
so we may assume in addition that 
\begin{equation}\label{grad-assumption}
\partial^{p_0}_s P(0,0) \ = \ \partial_t P(0,0) \ = \ 0.
\end{equation}
This innocent looking reduction will be fundamental later on, allowing us to estimate certain oscillatory integrals efficiently.

For Theorem \ref{main-uniform}, we take $p_0 = 1$, 
$P$ a general real polynomial and ${\mathcal F}$ a general
finite set; our goal is to obtain $L^2({\mathbb H}^1)$ bounds, uniform with respect to ${\mathcal F}$ and 
$P$ lying in some subspace ${\mathcal V}_{\Delta}$ of real polynomials. For Theorem \ref{vw-analytic} we consider
general $p_0\ge 1$ and real-analytic $P$ near $(0,0)$, but we insist that the dyadic rectangles
${\mathcal R}_I$ associated to $I \in {\mathcal F}$ all lie in some small fixed neighbourhood (depending on $P$)
of the origin $(0,0)$; no uniformity in $P$ is sought in our $L^2$ bounds for the corresponding operators.

In analysing $T_{P, {\mathcal F}}$ we take an oscillatory integral approach.
Viewing $T_{P,{\mathcal F}} f = L \ast_{\mathbb{H}^1} f$ as a Heisenberg convolution operator, one
can deduce via the group Fourier transform
on ${\mathbb H}^1$, that
\[
\|T_{P, {\mathcal F}}\|_{L^2({\mathbb H}^1) \to L^2({\mathbb H}^1)} \ \sim \ \ \sup_{\lambda \in {\mathbb R}}
\|S_{P, {\mathcal F}}^{\lambda}\|_{L^2({\mathbb R}) \to L^2({\mathbb R})}
\]
where
\[
S_{P, {\mathcal F}}^{\lambda} g (x) \ = \ \sum_{I \in {\mathcal F}} \int_{\mathbb R} m_I (\lambda, y,t) g(t) \, dt \ =: \
\sum_{I \in {\mathcal F}} S^{\lambda}_{P,I} g(y)
\]
and
\[
m_I(\lambda, y, t) \ = \ \int_{\mathbb R} e^{2\pi i \lambda ( (y+t) s^{p_0} + P(y-t, s))} \,  \phi_I^{(I)}(s,y-t) \, ds.
\]
See \cite{S} for an expression for the Fourier transform on $\mathbb{H}^1$.

\begin{remark}
Here we must caution the reader that the above reduction to a multiplier question on $\widehat{\mathbb{H}^1}$ does not come for free. Indeed, $L$ above is a distribution and there is no a priori reason for it to have a well-behaved group Fourier transform. However, with a little care one can verify that the above reduction is indeed justified. For details, see for example \cite{Kim} where an analogous one-parameter singular Radon transform is considered.
\end{remark}

If $P(y-t, s) = \sum_{p,q \ge 0} c_{p,q} (y-t)^q s^p$, then since $P(0,0) = 0$ and \eqref{grad-assumption} holds, we can write
$P(y-t, s) = \varphi(s) + \psi_0(y-t) + \sum_{p\ge 1} \psi_p(y-t) s^p$ where
\[
\varphi(s) = \sum_{p\ge 1} c_{p,0} s^p, \ \ \psi_0(y-t) = \sum_{q\ge 2} c_{0,q} (y-t)^q \ \ {\rm and}
\ \ \psi_p(y-t) = \sum_{q\ge 1} c_{p,q} (y-t)^q
\]
so that $\psi_p(0) = 0$ for all $p\ge 0$ (and $\psi_0'(0) = 0$). Importantly we have $c_{p_0,0} = 0$.

We can write the phase $(y+t) s^{p_0} + P(y-t, s)$ of $m_I$ as
$2 y s^{p_0} + {\tilde P}(y-t, s)$ where the difference between $P(y-t, s)$ and ${\tilde P}(y-t,s)$
is that the coefficient $c_{p_0, 1}$ in $\psi_{p_0}(y-t)$ is changed to $c_{p_0, 1} - 1$. This change does
not affect $P_{p_0}$ and so in the proofs of either Theorems \ref{main-uniform} or \ref{vw-analytic}
we may assume, without loss of generality, that 
\begin{equation}\label{mI}
m_I(\lambda, y, t) \ = \ \int_{\mathbb R} e^{2\pi i \lambda ( 2 y s^{p_0} + P(y-t, s))} \,  \phi_I^{(I)}(s,y-t) \, ds.
\end{equation}

Clearly bounds on the oscillatory integral $\sum_{I \in {\mathcal F}} m_I(\lambda, y, t)$
will play a central r\^{o}le in our analysis. General estimates for oscillatory integrals will be detailed in the
next section but for now we highlight a couple generalisations of an important, well-known oscillatory integral
bound due to Stein and Wainger \cite{SW} which states that for any real polynomial $Q \in {\mathbb R}[X]$,
we have
\[
\Bigl| \int_{a\le |s| \le b} e^{2 \pi i Q(s)} \, \frac{ds}{s} \Bigr| \ \le \ C_d
\]
where $C_d$ depends only on the degree $d$ of $Q$ and otherwise independent of the coefficients 
of $Q$ as well as $a$ and $b$. A proof from a modern perspective is given in \cite{S} and in fact gives the stronger
bound
\[
\sum_{k \in S} |n_k(Q)| \ \le \ C_d \ \ {\rm where} \ \ n_k(Q) \ = \ \int_{|s|\sim 2^k} e^{2\pi i Q(s)} \, \frac{ds}{s}
\]
and $S$ is any set of integers and $C_d$ can be taken to be independent of $S$. In our context, we need to show
that for any subset ${\mathcal F}' \subseteq {\mathcal F}$,
\begin{equation}\label{sw-generalised}
\sum_{I \in {\mathcal F}'} | m_I(\lambda, y, t) | \ \le \ C \, \frac{1}{|y-t|}
\end{equation}
holds when either (i) $P$ is a general real polynomial and $C = C_d$ depends only on the degree $d$
of $P$ (and in particular does not depend on the subset ${\mathcal F}' \subset {\mathbb Z}^2$,
$\lambda$, $y$, $t$ and the coefficients of $P$)
or (ii) $P$ is real-analytic near $(0,0)$ and ${\mathcal F}$ indexes dyadic rectangles ${\mathcal R}_I$
supported near the origin; that is,
the pairs $I = (j,k)$ range over integers $j\le -J$ and $k\le -K$ where $J$ and $K$ are
large, fixed positive integers depending on our real-analytic function $P$. In this case, the constant
$C$ is allowed to depend on $P$ and in particular it will depend on the truncation parameters $J, K$
but it does not depend on $\lambda, y, t$ or the cardinality of ${\mathcal F}'$. 

In the next section we will establish the estimate \eqref{sw-generalised} in both cases.

\subsection{Hilbert integral reduction} Choose $\chi \in C^{\infty}_0({\mathbb R})$ supported in
$\{|y| \sim 1\}$ and such that
if $\chi_m (y) := \chi(2^{-m}y)$, we have $\sum_{m\in {\mathbb Z}} \chi_m(y) \equiv 1$ for $y\neq 0$. 
We decompose 
\[
S^{\lambda}_{P,{\mathcal F}} g(y) \ = \ \sum_{m\in{\mathbb Z}} \sum_{I \in {\mathcal F}} \chi_m(y) 
S^{\lambda}_{P,I} g(y)  \ = \ S^1 g(y) + S^2 g(y)
\]
where
\[
S^1 g(y) \ := \ \sum_{(m,I) \in {\mathcal L}^1} \chi_m(y) S^{\lambda}_{P,I} g(y) \ \ {\rm and} \ \
S^2 g(y) \ := \  
 \sum_{(m,I) \in {\mathcal L}^2} \chi_m(y) S^{\lambda}_{P,I} g(y)
\]
and 
\[
{\mathcal L}^1  \ := \ \{(m,I) \in {\mathbb Z}\times {\mathcal F}: I = (j,k) \ {\rm satisfies} \ m \le k + C_0 \}
\]
for some large, fixed $C_0 > 0$. The set ${\mathcal L}^2$ is defined similarly but with the condition $k \le m - C_0$. The significance of this is that when $(m,I) \in \mathcal{L}^2$ we have $2^k \sim |y-t| \ll |y|\sim 2^m$.\\
Hence
\[
|S^1 g(y)| \ \le \ \int \bigl[ \sum_{(m,I) \in {\mathcal L}^1} |\chi_m(y) m_I(\lambda, y, t)| \bigr] \, |g(t)| \, dt 
\]
where the sum over $(m,I) \in {\mathcal L}^1$
is supported in $\{ (y,t): \delta |y| \le |y-t| \}$ for some small $\delta>0$, depending on our choice of
$C_0$. Using \eqref{sw-generalised}, we have
\[
|S^1 g(y)| \ \le \ C \, \int_{\delta |y| \le |y-t| } \frac{1}{|y-t|} \, |g(t)| dt \ =: \ \int K(y,t) |g(t)| \, dt.
\]
The integral operator with kernel $K$ is of {\it Hilbert integral} type (the kernel is homogeneous of degree
$-1$ and $K(1,t) |t|^{-1/2}$ is integrable over ${\mathbb R}$) and hence $S^1$ is uniformly bounded
on $L^2({\mathbb R})$ (uniform in $\lambda$, ${\mathcal F}$ and the coefficients of $P$ in the polynomial case).
See \cite{S1}, page 271.

For $S^2$, write $S^2 g (y) = \sum_m S^2_m g(y)$ where
\[
S^2_m g(y) \ = \ \sum_{I: (m,I) \in {\mathcal L}^2} \chi_m(y) S^{\lambda}_{P,I} g(y) \ =: \ \int K_m (y,t) g(t) \, dt.
\]
For $|y| \sim 2^m$, we have $|t| = |t-y + y| \sim |y| \sim 2^m$ if $|y-t| \sim 2^k$
and $k\le m - C_0$. Hence ${\rm supp}(K_m) \subset \{(y,t): |y|, |t| \sim 2^m \}$ and so
\begin{equation}\label{Sm}
\|S^2 \|_{L^2 \to L^2} \ = \ \|\sum_m S^2_m \|_{L^2 \to L^2} \ \sim \ \sup_m \, \|S^2_m \|_{L^2 \to L^2}
\end{equation}
by (almost) orthogonality. Therefore the proofs of both Theorem \ref{main-uniform}
and Theorem \ref{vw-analytic} reduce to understanding when the operators $S^2_m$
are uniformly bounded on $L^2$.

\section{Oscillatory integral estimates}\label{oscillatory}

Many oscillatory estimates rely on van der Corput's lemma which we now state.

{\bf van der Corput's Lemma} {\it For any $k\ge 2$, there exists a constant $C_k$ such that
\[
\Bigl| \int_a^b e^{2\pi i \lambda \phi(s) } ds \Bigr| \ \le \ C_k |\lambda|^{-1/k} 
\]
holds for any real-valued $\phi \in C^k[a,b]$ such that $|\phi^{(k)}(s)| \ge 1$ for $s \in [a, b]$.
The result holds for $k=1$ if in addition we assume that $\phi'$ is monotone on $[a,b]$.}

For a proof, see \cite{S}. 
Let $Q(s) = \lambda [ 2 y s^{p_0} + P(y-t,s)]$ be the phase appearing in each 
\[m_I(\lambda, y, t) \ = \
 \int_{\mathbb R} e^{2\pi i Q(s)} 2^{-j-k} \phi_I(2^{-j} s, 2^{-k}(y-t)) \, ds \ =: \ 2^{-k} I_j
\]
where
\begin{equation}\label{Ij}
I_j \ = \ I_{j,\lambda, y, t, k} \ = \ \int_{\mathbb R} e^{2\pi i Q(s)} \, 2^{-j} \Phi(2^{-j} s)  \, ds 
\end{equation}
is supported in $\{(y,t,k) \;:\; |y-t| \sim 2^k \}$. Here the $\Phi = \Phi_{y,t,I, k}$
have bounded $C^{\ell}$ norms, uniformly in the parameters $y,t, I$ and $k$. Hence
$|m_I(\lambda, y, t)| \le |I_j| \ 2^{-k} \, \chi_{|y-t|\sim 2^k}$ and
so to bound $\sum_I |m_I(\lambda, y, t)|$, it suffices to fix $k$ and obtain uniform bounds for
the sums $\sum_j |I_j|$ over $j$. To do this, 
we will use
van der Corput's lemma.

Our first application is a proof of \eqref{sw-generalised}.

\subsection{Proof of the generalised Stein-Wainger bound \eqref{sw-generalised}} 

Let $Q(s) = \lambda [ 2 y s^{p_0} + P(y-t,s)]$ be the phase appearing in each $m_I$ and
for each $k\in {\mathbb Z}$, set ${\mathcal F}_k' = \{j \in {\mathbb Z} : I = (j,k) \in {\mathcal F}'\}$.
It suffices to show that for every $k\in {\mathbb Z}$,
\begin{equation}\label{k}
\sum_{I\;:\; j\in {\mathcal F}_k'} |m_I(\lambda, y, t)| \ \lesssim \  \, 2^{-k} \chi_{|y-t| \sim 2^k }
\end{equation}
since \eqref{sw-generalised} follows by summing these estimates over $k\in {\mathbb Z}$. As observed above,
this is equivalent to showing
\[
\sum_{j \in {\mathcal F}_k'} |I_j| \ \lesssim \ 1
\]
where
\[
I_j \ := \ \int_{\mathbb R} e^{2\pi i Q(s)} \, 2^{-j} \Phi(2^{-j} s)  \, ds \ = \ \int_{\mathbb R} 
e^{2\pi i Q_j(s)} \, \Phi(s) \, ds
\]
and $Q_j(s) = Q(2^j s)$.

We start with
the case when $P$ is a polynomial where we seek bounds which are uniform in the coefficients of $P$,
the subset ${\mathcal F}' \subseteq {\mathcal F}$, 
and the parameters $\lambda, y$ and $t$. For the case when $P$ is real-analytic at $(0,0)$,
we will reduce the estimate \eqref{sw-generalised} to the polynomial case.

In the polynomial case, our phase $Q(s) = \sum_{p\ge 1} e_p s^p$
is a polynomial (without loss of generality we may suppose that $Q$ has no constant term) and
hence $Q_j(s) = \sum_{p\ge 1} e_p 2^{p j} s^p$. 
A simple equivalence of norms argument shows that there exists a $c_d >0$, depending only on
the degree $d$ of $Q$, such that for all $j$ there exists $\ell_j$ with $1 \leq \ell_j \leq d$ for which
$|Q_j^{(\ell_j)}(s)| \ge c_d \sum_{p\ge 1} |e_p| 2^{p j}$ holds on the support of $\Phi$. An application of
van der Corput's lemma now shows that $|I_j| \le C_d (\Lambda_j)^{-1/\ell_j}$ where
$\Lambda_j = \sum_{p\ge 1} |e_p| 2^{p j}$. Using \eqref{product-cancellation},   
\[
\int_{\mathbb R} \Phi(s) \, ds \ = \ \int_{\mathbb R} \phi_I (s, 2^{-k}(y-t)) \, ds \ = \ 0,
\]
one also has
\[
|I_j| = \Big|\int_{\mathbb R} [e^{2\pi i Q_j(s)} - 1] \Phi(s) ds \Big| \ \lesssim \ \Lambda_j \ \ {\rm and \ so} \ \
|I_j| \lesssim \min(\Lambda_j, \Lambda_j^{-1/d}),
\]
which allows us to sum in $j$ to see $\sum_j |I_j| \lesssim_d 1$, as desired.

Next we consider the real-analytic case so that the pairs $I=(j,k)$ in ${\mathcal F}$
range over integers $j\le -J$ and $k\le -K$ where $J$ and $K$ are
large, fixed positive integers depending on our real-analytic function $P$. 
In this case, as said above, we will reduce matters to the polynomial case. 

Recall our notation where we write $P(y-t,s) = \phi(s) + \sum_{p\ge 0} \psi_p(y-t) s^p$  and
$\phi(s) = \sum_{p\ge 1} c_{p,0} s^p$. For $|y-t| \sim 2^k \ll 1$,  we have $|\psi_p(y-t)| \lesssim_P 2^{k}$.
Write
\[
|m_I(\lambda, y, t)| = \Bigl| \int_{\mathbb R} e^{2\pi i \lambda [ 2 y s^{p_0} + \phi(s) + \Psi_{y,t}(s)]} 
\, 2^{-j} \Phi(2^{-j} s) \, ds\Bigr|
\]
where $\Psi_{y,t}(s) := \sum_{p\ge 1} \psi_p(y-t) s^p$. Hence
$|\partial^{(p)}_s \Psi_{y,t}(s)| \lesssim_{p, P} \, 2^{k}$ for every $p\ge 0$ and in particular,
$|\partial^{(p)}_s \Psi_{y,t}(s)| \ll_{p, P} \, 1$.

First we consider the case that there exists a $p_1 > p_0$ such that $c_{p_1, 0} \neq 0$. 
Hence
$|\phi^{(p_1)}(s)| \gtrsim_P 1$ for $|s| \ll 1$ and so 
\[
|\partial_s^{(p_1)} [ 2 y s^{p_0} + \phi(s) + \Psi_{y,t}(s) ] | \ \gtrsim \ 1.
\]
This puts us in a position to apply van der Corput's lemma, which together with a simple integration by
parts argument allows us to conclude $|I_j| \lesssim 2^{-j} |\lambda|^{-1/p_1}$ 
and so $\sum_{j\in S_{\lambda}} |I_j| \lesssim 1$ where
$S_{\lambda} = \{ j: 2^j \ge |\lambda|^{-1/p_1} \}$. For $j\notin S_{\lambda}$, we compare
the integral $I_j$ to the integral 
\[
II_j \ := \ \int_{\mathbb R} e^{2\pi i \lambda [ 2 y s^{p_0} + {\tilde \phi}(s) + {\tilde \Psi}_{y,t}(s)]} 
\, 2^{-j} \Phi(2^{-j} s) \, ds
\]
where ${\tilde \phi}(s) = \sum_{p=1}^{p_1 - 1} c_{p,0} s^p$ and 
${\tilde \Psi}_{y,t}(s) = \sum_{p=1}^{p_1 -1} \psi_p(y-t) s^p$. Note that the difference of the phases
in $I_j$ and $II_j$ is at most $C |\lambda s^{p_1}|$ and so
\[
|I_j - II_j| \ \lesssim \ |\lambda| 2^{p_1 j},
\]
implying $\sum_{j\notin S_{\lambda}} |I_j - II_j| \lesssim 1$. We can appeal to our analysis
of \eqref{sw-generalised} when the phase is
polynomial to conclude $\sum_{j\notin S_{\lambda}} |II_j| \lesssim 1$ and hence
\eqref{k} holds in this case.

Finally we consider the case $\phi(s) = \sum_{p< p_0} c_{p,0} s^p$; that is, there is no $p_1 > p_0$ such
that $c_{p_1, 0} \neq 0$ (remember $c_{p_0,0} = 0$ by \eqref{grad-assumption}). In this case we may suppose that there is a $p_1 > p_0$ such that
$|\psi_{p_1}(y-t)| \sim |y-t|^{\ell_{*}}$ for some $\ell_{*} \ge 1$ and $\psi_p^{(\ell)}(0) = 0$
for all $p\ge p_1$ and all $\ell < \ell_{*}$. Otherwise $\psi_p \equiv 0$ for all $p > p_0$ and we
are back in the polynomial case. In particular
\[
\partial^{(p_1)}_s [2 y s^{p_0} + \phi(s) + \Psi_{y,t}(s) ]  = c (y-t)^{\ell_{*}} + O( (y-t)^{\ell_{*}} s)
\]
and therefore $|\partial^{(p_1)}_s [2 y s^{p_0} + \phi(s) + \Psi_{y,t}(s) ]| \gtrsim 2^{\ell_{*} k}$
for $|s| \ll 1$ and $|y-t| \sim 2^k$. Hence by van der Corput's lemma, 
$|I_j| \lesssim 2^{-j} [|\lambda| 2^{\ell_{*} k}]^{-1/p_1}$ implying that
$\sum_{j\in S_{\lambda}'} |I_j| \lesssim 1$ where 
$S_{\lambda}' := \{j: 2^j \ge (|\lambda| 2^{\ell_{*} k})^{-1/p_1}\}$. For $j \notin S_{\lambda}'$,
we compare
the integral $I_j$ to the integral 
\[
III_j \ := \ \int_{\mathbb R} e^{2\pi i \lambda [ 2 y s^{p_0} + \phi(s) + {\tilde \Psi}_{y,t}(s)]} 
\, 2^{-j} \Phi(2^{-j} s) \, ds
\]
where again
${\tilde \Psi}_{y,t}(s) = \sum_{p=1}^{p_1 -1} \psi_p(y-t) s^p$. Note that the difference of the phases
in $I_j$ and $III_j$ is at most $C |\lambda 2^{\ell_{*} k} s^{p_1}|$ and so
\[
|I_j - III_j| \ \lesssim \ |\lambda| 2^{\ell_{*} k} 2^{p_1 j},
\]
implying $\sum_{j\notin S_{\lambda}'} |I_j - III_j| \lesssim 1$.
Once again we can appeal to our analysis
of \eqref{sw-generalised} when the phase is
polynomial to conclude $\sum_{j\notin S_{\lambda}'} |III_j| \lesssim 1$ and hence
\eqref{k} holds in this case as well. This completes the proof of \eqref{sw-generalised} in all cases.

\subsection{Another useful bound for oscillatory integrals}\label{} A nontrivial application of van der Corput's lemma
gives the following useful uniform bound for oscillatory integrals with polynomial phases.

\begin{proposition}\label{nw} For any $Q(s) = \sum_{j=1}^d h_j s^j \in {\mathbb R}[X]$ and $1\le j \le d$,
we have
\[
\Bigl| \int_{B/2}^B e^{2\pi i Q(s) } \, ds \Bigr| \ \le \ C_d |h_j|^{-1/d} B^{1 - j/d}.
\]
\end{proposition}

This is a simple variant of Theorem 3.1 in \cite{NW}. We have the following immediate consequence for our
multipliers $m_I(\lambda, y, t) = 2^{-k} I_j$ when the phase
$Q(s) = \lambda ( 2 y s^{p_0} + P(y-t, s) ) =  \sum_{p \le d} h_p s^p$ is a polynomial. We have for
every $1\le p \le d$,
\begin{equation}\label{Ij-poly}
|I_j|  \lesssim_d  \bigl[|h_p| 2^{p j}\bigr]^{- 1/d}
\end{equation}
where we recall the definition of $I_j$ in \eqref{Ij}.
We will use this estimate in the proof
of Theorem \ref{main-uniform} where $P$ is a polynomial. For Theorem \ref{vw-analytic},
when $P$ is assumed to be real-analytic near $(0,0)$, we will need the following two variants of \eqref{Ij-poly}.

Consider again the phase 
\[
Q(s) = \lambda ( 2 y s^{p_0} + P(y-t,s) ) =
 \lambda \Bigl[ (2 y s^{p_0} + \sum_{p\ge 1} c_{p,0} s^p +  \sum_{p \ge 0}  \psi_p(y-t) s^p \Bigr]
\]
in $m_I(\lambda, y, t) = 2^{-k} I_j$. The coefficient of $s^{p_0}$ is 
$h_{p_0} = \lambda (2 y + \psi_{p_0}(y-t))$, again, since the coefficient $c_{p_0,0} = 0$ as per \eqref{grad-assumption}. 
This is important since it allows the size of $h_{p_0}$ to be determined in the setting of $S^2_m$ where
matters have been reduced (see \eqref{Sm}).

We consider $I_j = 2^k m_I(\lambda, y, t)$ where the
the pair $(m,I) \in {\mathcal L}^2$ arises in the definition of $S^2_m$. Hence 
$|y-t| \sim 2^k \ll 2^m \sim |y|$ and so
$|h_{p_0}| \sim |\lambda| 2^m$
since $\psi_{p_0}(y-t) = O_P(2^k)$. In this case,  we have
\begin{equation}\label{Ij-analytic-0}
|I_j| \ \lesssim_{P} \ \bigl[|h_{p_0}| 2^{p_0 j}\bigr]^{-\epsilon} \ \sim \ 
\bigl[|\lambda| 2^m 2^{p_0 j}\bigr]^{-\epsilon}
\end{equation}
for some $\epsilon> 0$.

Next we consider an estimate with respect to the coefficient $h_p = \lambda (c_{p,0} + \psi_p(y-t))$
of $s^p$ in the phase $Q(s)$ for other values of $p$. In our arguments, this case will only arise
in the simpler situation when the phase $Q$  is truncated to either
\[
 \lambda \Bigl[ (2 y s^{p_0} + \sum_{p\ge 1} c_{p,0} s^p +  \sum_{1\le p < p_0}  \psi_p(y-t) s^p \Bigr]
\]
or
\[
 \lambda \Bigl[ (\sum_{p\ge 1} c_{p,0} s^p +  \sum_{1\le p < p_0}  \psi_p(y-t) s^p \Bigr]
\]
which is still not quite the case of a polynomial.
For any $1\le p < p_0$ with $\psi_p \not\equiv 0$, we have for some $\ell_p \ge 1$, 
$|\psi_p(y-t)| \sim 2^{\ell_	p k}$ when $|y-t| \sim 2^k$. Hence the coefficient of
$h_p = \lambda( c_{p,0} + \psi_p(y-t))$ satisfies $|h_p| \gtrsim 2^{\ell_p k}$ since
$2^{\ell_p k} \ll 1$. In this situation, we have
\begin{equation}\label{Ij-analytic-p}
|I_j| \ \lesssim_{P,p} \ \bigl[|h_{p}| 2^{p j}\bigr]^{-\epsilon} \ \lesssim \ 
\bigl[|\lambda| 2^{\ell_p p} 2^{p j}\bigr]^{-\epsilon}
\end{equation}
for some $\epsilon> 0$.

The proof of \eqref{Ij-analytic-0} is fairly simple and we present this case now. The proof of \eqref{Ij-analytic-p} is an elaboration on a proof of Proposition \ref{nw} and we have decided to give the proof in an appendix to the paper.

To prove \eqref{Ij-analytic-0}
we begin as in the real-analytic case for \eqref{k} by initially assuming
there exists a $p_1 > p_0$ such that $c_{p_1, 0} \neq 0$. 
Hence
$|\phi^{(p_0)}(s)| \sim |s|^{p_1 - p_0}$ for $|s| \ll 1$ and so if $|s| \ll |y|^{1/(p_1 - p_0)}$
or $|y|^{1/(p_1 - p_0)} \ll |s|$  (that is, $2^j \not\sim 2^{m/(p_1 - p_0)}$), we see that 
\[
|Q^{(p_0)}(s)| = |\lambda ( 2 p_0 ! \, y + \phi^{(p_0)}(s) + O(2^k) )| \ \gtrsim \
|\lambda| \bigl[ |y|  - C 2^k \bigr] \ \gtrsim |\lambda| |y| \ \sim \ |h_{p_0}|
\]
since $2^k \ll 2^m \sim |y|$. Hence by van der Corput's lemma, we have 
$|I_j| \lesssim_P 2^{-j} |h_{p_0}|^{-1/p_0}$ implying \eqref{Ij-analytic-0} with
$\epsilon = 1/p_0$.

When $2^j \sim  2^{m/(p_1 - p_0)}$, we consider the $p_1$th derivative of $Q$:
note that $|\phi^{(p_1)}(s)| \sim_P 1$ for $|s| \ll 1$. Therefore we have
\[
|Q^{(p_1)}(s)| = |\lambda ( \phi^{(p_1)}(s) + O(2^k) )| \ \gtrsim \
|\lambda| 
\]
since $2^k \ll 1$. Hence van der Corput's lemma implies 
\[
|I_j| \ \lesssim_P \ 2^{-j} |\lambda|^{-1/p_1} \ = \ (|\lambda| 2^{j(p_1-p_0)} 2^{p_0 j})^{-1/p_1} \ \sim \
( |h_{p_0}| 2^{p_0 j})^{-1/p_1},
\]
implying \eqref{Ij-analytic-0} with $\epsilon = 1/p_1$.

Finally we consider the case that for all $p_1 > p_0$ we have $c_{p_1, 0} = 0$ in which
case $\phi^{(p_0)}(s) \equiv 0$ since we also have $c_{p_0} = 0$. Therefore as before,
\[
|Q^{(p_0)}(s)| = |\lambda ( 2 p_0 ! \, y  + O(2^k) )| \ \gtrsim \
|\lambda| \bigl[ |y|  - C 2^k \bigr] \ \gtrsim |\lambda| |y| \ \sim \ |h_{p_0}|
\]
since $2^k \ll 2^m \sim |y|$. Hence by van der Corput's lemma, we have 
$|I_j| \lesssim_P (|h_{p_0}| 2^{p_0 j})^{-1/p_0}$ implying \eqref{Ij-analytic-0} with
$\epsilon = 1/p_0$. This completes the proof of \eqref{Ij-analytic-0} in all cases.

\section{The proof of Theorems \ref{main-uniform} and \ref{vw-analytic} -- the main steps}\label{vw}
In both Theorems \ref{main-uniform} and \ref{vw-analytic}, we need to establish uniform (in $m$) $L^2$
bounds for the operators
\[
S^2_m g(y) =  \sum_{I: (m,I) \in {\mathcal L}^2} \chi_m(y) S^{\lambda}_{P,I} g(y) \ \ {\rm where} \ \
S^{\lambda}_{P,I} g (y) = \int_{\mathbb R} m_I (\lambda, y,t) g(t) \, dt
\]
and
\[
m_I(\lambda, y, t)  \ =  \ \int_{\mathbb R} e^{2\pi i \lambda ( 2 y s^{p_0} + P(y-t, s))} \,  \phi_I^{(I)}(s,y-t) \, ds.
\]
See \eqref{Sm}. Here 
\[
{\mathcal L}^2 \ = \ \{(m,I) \in {\mathbb Z}\times {\mathcal F}: I = (j,k) \ {\rm satisfies} \ k \le m - C_0  \}
\]
for some large, fixed $C_0 > 0$. Recall that we write
\[
P(y-t, s) = \phi(s) + \sum_{p\ge 0} \psi_p(y-t) s^p  =  \phi(s) + \psi_0(y-t) + \sum_{p\in {\mathcal P}} \psi_p(y-t) s^p
\]
where each $\psi_p(0) = 0$ and ${\mathcal P} := \{p\ge 1: \psi_p \not\equiv 0\}$.

The plan of the proof is to use the oscillatory integral estimates discussed in Section \ref{oscillatory} to bound the errors introduced when removing certain terms from the phase of $m_I$. We will keep removing terms from the phase whenever possible until we have reduced matters to (euclidean convolution) operators that are well-known already. These will be either (variable kernel) oscillatory singular integral operators \`{a} la Ricci-Stein \cite{RS1} or the singular Radon transforms mentioned in the statements of Theorems \ref{main-uniform} and \ref{vw-analytic}.


\subsection{The exceptional set ${\mathcal E}$}
For both theorems, we will need to avoid an exceptional set ${\mathcal E}$ of bad values of $k$ which we
will make more and more explicit as we proceed. 
For Theorem \ref{main-uniform}, the cardinality $\# {\mathcal E} \lesssim_d 1$
is bounded uniformly in ${\mathcal F}$ and the coefficients of $P$. For Theorem \ref{vw-analytic},
the cardinality $\# {\mathcal E} \lesssim_P 1$ depends on $P$ (and hence on the truncation parameters $J,K$)
but is otherwise independent of $\#{\mathcal F}$.

We split 
\[
S^2_m g(y) \ = \ \sum_{I \in {\mathcal F}^{0,m}}  \chi_m(y) S^{\lambda}_{P,I} g(y)  +  
\sum_{I \in {\mathcal F}^{1,m}} \chi_m(y) S^{\lambda}_{P,I} g(y) \ := \ S^{2,0}_m g(y) + S^{2,1}_m g(y)
\]
where ${\mathcal F}^{0,m} = \{ I = (j,k) : (m,I) \in {\mathcal L}^2, k \notin {\mathcal E} \}$ and
${\mathcal F}^{1,m}$ involves the bad values $k \in {\mathcal E}$.  We use \eqref{k} with 
${\mathcal F}' = {\mathcal F}^{1,m}$
to bound
\[
|S^{2,1}_m g(y)| \le \sum_{k\in {\mathcal E}} \int \bigl[\sum_{j\in{\mathcal F}^{1,m}_k} 
|m_I(\lambda, y, t)| \bigr] \, |g(y)| dy \ \lesssim \ \sum_{k\in {\mathcal E}}  2^{-k} \int_{|y-t| \sim 2^k} 
|g(y)| \, dy
\]
and so $\|S^{2,1}_m\|_{L^2 \to L^2} \lesssim \# {\mathcal E} \lesssim 1$, leaving us with
$S^{2,0}_m$ which avoids the bad values $k \in {\mathcal E}$. 

To ease the notation, we rewrite
$S^{2,0}_m$ as $S^2_m$ with the understanding that the sum defining $S^2_m$
is taken over $I = (j,k) \in {\mathcal F}^{0,m}$ and so $k \notin {\mathcal E}$. 

For each term $\psi_{p_{*}}(y-t)$ with $p_{*} \in {\mathcal P}$ arising in the phase of $m_I$,
our strategy is to reduce the analysis of $S_m^2$ to 
$S_m^{2,*} = \sum_{I\in {\mathcal F}^{0,m}} \chi_m(y) S^{\lambda, *}_{P,I}$ where 
\[
S^{\lambda,*}_{P,I} g (y) \ = \ \int_{\mathbb R} m_I^{*} (\lambda, y,t) \   g(t) \, dt
\]
and
\begin{equation}\label{m*}
m_I^{*}(\lambda, y, t)  \ =  \ \int_{\mathbb R} 
e^{2\pi i \lambda ( 2 y s^{p_0} + \phi(s) + \sum_{p\not = p_{*}} \psi_p(y-t) s^p)} \,  \phi_I^{(I)}(s,y-t) \, ds;
\end{equation}
that is, we plan to remove the term $\psi_{p_\ast}(y-t)s^{p_\ast}$ from the phase of $m_I$.

Our estimates are naturally expressed in terms of certain key quantities associated to the size of those
$\psi_{p_{*}}(y-t)$ with $p_{*} \in {\mathcal P}$. For Theorem \ref{vw-analytic}, when
$P$ is assumed to be real-analytic near $(0,0)$, we can find
an $\ell_{*}\ge 1$ such that $|\psi_{p_{*}}(y-t)| \sim c_{*} 2^{\ell_{*} k}$ when $|y-t| \sim 2^k \ll 1$.
This simply follows from the fact that $\psi_{p_{*}}(0) = 0$ and $\psi_{p_{*}} \not\equiv 0$. 
For Theorem \ref{main-uniform} the $\psi_p(y-t)$ are general polynomials and $|y-t| \sim 2^k$ can be
of any size ($k \in {\mathbb Z}$ can take any value). Here we will appeal to a result in \cite{CRW} which shows
that outwith finitely many values of $k$ (depending only on the degree of $P$), there exists
an  $\ell_{*}\ge 1$ such that indeed $|\psi_{p_{*}}(y-t)| \sim c_{*} 2^{\ell_{*} k}$ when $|y-t| \sim 2^k$. 

Given a nonzero polynomial $Q \in {\mathbb R}[X]$, a basic result in \cite{CRW} gives us a decomposition
${\mathbb R} = S \cup G$ where $S = \cup I$ can be written as a disjoint union of $O(1)$ (with constant
only depending on the degree of $Q$) intervals such that on each $I$, 
$|Q(t)| \sim c_I |t|^{\ell_I}$ for some $\ell_I \in {\mathbb N}$. Furthermore if $Q(0) = 0$, then
$\ell_I \ge 1$ for all $I$. Finally each interval comprising $G = {\mathbb R}\setminus S$ is a dyadic
interval of the form $[A, C A]$ where $C \lesssim 1$. 

As above, we write our polynomial 
$P$ as $P(s,t) = \phi(s) + \sum_{p\ge 0} \psi_p (t) s^p$ where each $\psi_p \in {\mathbb R}[X]$
satisfies $\psi_p(0) = 0$. We apply the decomposition in \cite{CRW} to each $\psi_p $ with $p\in {\mathcal P}$
(so that $\psi_p \not\equiv 0$) to conclude
that there is an exceptional set ${\mathcal B}$ of $O_d(1)$ values of $k$ where
${\mathbb Z}\setminus{\mathcal B} = \cup_{n=1}^N S_n$ decomposes
into $O_d(1)$ sets such for each $p\in{\mathcal P}$ and $n$,
there is an $\ell_p = \ell_{p,n}\ge 1$ and $c_p = c_{p,n} > 0$ with the property that
\begin{equation}\label{dyadic-size}
|\psi_p(y-t)| \ \sim \ c_{p} \, 2^{\ell_{p} k} \ \ {\rm whenever} \ \ |y-t| \sim 2^k \ \ {\rm and} \ \ k \in S_n.
\end{equation}
We incorporate the set ${\mathcal B}$ into ${\mathcal E}$ so that $I = (j,k) \in {\mathcal F}^{0,m}$
implies $k \in S_n$ for some $n$ and \eqref{dyadic-size} holds for every $\psi_p$ with $p\in {\mathcal P}$.

\subsection{Key quantities and the first step}\label{key-step}
The key quantities
${\mathcal A}_{p_{*}}(k) = {\mathcal A}_{p_{*},\lambda, m}(k)$ are defined as
\[
{\mathcal A}_{p_{*}}(k) \ := \ \frac{|\lambda| c_{*} 2^{\ell_{*} k}}{(|\lambda| 2^m)^{p_{*}/p_0}}
\]
where, in the case of Theorem \ref{main-uniform}, $c_{*} = c_{p_{*}}$ and $\ell_{*} = \ell_{p_{*}}$
appear in \eqref{dyadic-size}. 
One important estimate where these quantities arise occurs in the following bound for 
the differences $D_k := \sum_{j  \in {\mathcal F}^{0,m}_k} [m_I - m_I^{*}]$ (which avoids
the exceptional values of $k \in {\mathcal E}$),
\begin{equation}\label{Dk}
|D_k| \ \lesssim \ {\mathcal A}_{p_{*}}(k)^{\epsilon_{*}} \ \,  2^{-k} \chi_{|y-t| \sim 2^k}
\end{equation}
for some $\epsilon_{*}>0$.


For Theorem \ref{main-uniform}, the implicit constant in the estimate \eqref{Dk} will be uniform; it will depend only on
the degree of $P$ and can be taken to be independent of the coefficients of $P$ as well as the set
${\mathcal F}$. For Theorem \ref{vw-analytic} the implicit constant will depend on $P$. 

To prove
\eqref{Dk}, we split 
\[
D_k \ = \ \sum_{j\in J_1} [m_I - m_I^{*}] \ + \ \sum_{j\in J_2} [m_I - m_I^{*}] \ := \ D_k^1 + D_k^2
\]
where $J_1 \sqcup J_2 = \{j : I = (j,k) \in {\mathcal F}^{0,m}\}$ and
\[
J_1 \ := \ \{ j : I = (j,k) \in {\mathcal F}^{0,m} \ {\rm and} \ \  2^j \le (|\lambda| 2^m)^{-1/p_0} \, 
{\mathcal A}_{p_{*}}^{-\sigma}(k)  \}
\] 
for some $\sigma > 0$ to be chosen later.
For $j\in J_1$, we use that the difference in the phases of $m_I$ and $m_I^{*}$ is at most
$C |\lambda| c_{*} 2^{\ell_{*} k} 2^{p_{*} j}$ (the constant $C$ being absolute/uniform) to conclude that
$|D_k^1| \le$
\[
\sum_{j \in J_1} |m_I - m_I^{*}| \ \lesssim \  2^{-k} \, \chi_{|y-t| \sim 2^k} \ |\lambda| c_{*} 2^{\ell_{*} k} \,
\sum_{j \in J_1} 2^{p_{*} j} \ \lesssim \ {\mathcal A}_{p_{*}}(k)^{\epsilon_{*}} \ \,  2^{-k} \chi_{|y-t| \sim 2^k}
\]
where $\epsilon_{*} = 1 - \sigma p_{*} > 0$ and we have chosen $\sigma < 1/p_{*}$; this shows that \eqref{Dk} holds
for $D_k^1$. 

For $D_k^2$, we treat $m_I$ and $m_I^{*}$ separately, bounding $|D_k^2| \le \sum_{j\in J_2} |m_I| +
\sum_{j\in J_2} |m_I^{*}|$. Recall that 
\[
|m_I| \ = \ 2^{-k} \, 
\Bigl| \int_{\mathbb R} e^{2\pi i \lambda (2 y s^{p_0} + \phi(s) + \sum_{p\ge 1}\psi_p(y-t) s^p)} 2^{-j} 
\Phi(2^{-j} s) \, ds \Bigr|.
\]
We will apply \eqref{Ij-poly} and \eqref{Ij-analytic-0} to
\[
Q(s) =  \lambda \bigl[ (2 y + \psi_{p_0}(y-t)) s^{p_0} + \phi(s) +  \sum_{p\neq p_0} \psi_p(y-t) s^p \bigr]
\]
with respect to the coefficient $h_{p_0} := \lambda (2 y + \psi_{p_0}(y-t))$ of $s^{p_0}$. 
Very importantly, we have reduced (see \eqref{grad-assumption}) to the case where
the coefficient $c_{p_0,0}$ in $\phi(s) = \sum_{p\ge 1} c_{p,0} s^p$ is zero!

If $\psi_{p_0}(y-t) \equiv 0$, then $h_{p_0} = 2\lambda y$ and if $\psi_{p_0}(y-t) \not\equiv 0$, then
for $k \notin {\mathcal E}$, we have $|\psi_{p_0}(y-t)| \sim c_0  2^{\ell_0 k}$ for some $\ell_0\ge 1$ 
when $|y-t| \sim 2^k$.  Hence there are only $O(1)$ values of $k$
where the bound $|h_{p_0}| \sim |\lambda y| \sim |\lambda| 2^m$ does not hold. We add these values
to the expectional set ${\mathcal E}$. Hence \eqref{Ij-poly} and \eqref{Ij-analytic-0} imply 
\[
|m_I| \ \lesssim \ \bigl[|\lambda| 2^m 2^{p_0 j}\bigr]^{-\epsilon_0} \ 2^{-k} \chi_{|y-t| \sim 2^k}
\]
for some $\epsilon_0 > 0$. The same argument shows that $|m_I^{*}|$ satisfies this estimate as well. Summing
over $j\in J_2$ establishes \eqref{Dk} for $D_k^2$ and hence $D_k$.

\subsection{An interlude -- some analysis specific to Theorem \ref{vw-analytic}}
For Theorem \ref{vw-analytic} (in which case both $j, k\le 0$ for $I = (j,k) \in {\mathcal F}$),
we claim that 
when $p_{*} > p_0$, the above differences $D_k$ also satisfy 
\begin{equation}\label{Dk-neg}
|D_{k} | \ \lesssim_P \ {\mathcal A}_{p_{*}}(k)^{-\epsilon_{*}} \ \, 2^{-k} \chi_{|y-t|\sim 2^k}
\end{equation}
where 
$\epsilon_{*} = p_0/(p_{*} -p_0) >0$. This, together with \eqref{Dk}, will allow us to remove all terms $\psi_p(y-t) s^{p}$ with $p \geq p_0$ from the phase of $m_I$.

The proof of \eqref{Dk-neg} is straightforward. We again use that the difference in the phases
of $m_I$ and $m_I^{*}$ is at most
$C |\lambda| 2^{\ell_{*} k} 2^{p_{*} j}$ (the constant $C$ being absolute/uniform) to conclude that
$|D_k| \le$
\[
\sum_{j : I = (j,k) \in {\mathcal F}^{0,m}} |m_I - m_I^{*}| \lesssim_P   |\lambda|  2^{\ell_{*} k} 
\Bigl[\sum_{j \le 0} 2^{p_{*} j} \Bigr] \, 2^{-k}  \chi_{|y-t| \sim 2^k} \lesssim |\lambda|  2^{\ell_{*} k} 
\,  2^{-k} \chi_{|y-t| \sim 2^k}.
\]
However for $I = (j,k) \in {\mathcal F}^{0,m}$ we have $k\le 0$ and $k\le m$ and hence it can be verified that
\[
|\lambda| 2^{\ell_{*} k} \ \le \ |\lambda| 2^{k} \ \le \ 
\Bigl[\frac{(|\lambda| 2^m)^{p_{*}/p_0}}{|\lambda| 2^{\ell_{*} k}}\Bigr]^{p_0/(p_{*} - p_0)}.
\]
Therefore 
$|\lambda|  2^{\ell_{*} k} \le {\mathcal A}_{p_{*}}^{-\epsilon_{*}}(k)$
and so \eqref{Dk-neg} follows.

Note that when $k\le 0$ and $I = (j,k) \in {\mathcal F}^{0,m}$ (and so $k\ll m$), we have
\[
{\mathcal A}_{p_0} \ = \ c_0 \, 2^{\ell_0 k} 2^{-m} \ \le \ c_0 \, 2^k 2^{-m} \ \ll 1.
\]
Putting \eqref{Dk} and \eqref{Dk-neg} together, we see that in the situation of Theorem \ref{vw-analytic} 
and when $p_{*} \ge p_0$,
the differences satisfy
\[
|D_{k} | \ \lesssim_P \ \min({\mathcal A}_{p_{*}}(k), {\mathcal A}_{p_{*}}(k)^{-1})^{\epsilon_{*}} 
\ \, 2^{-k} \chi_{|y-t|\sim 2^k}
\]
for some $\epsilon_{*} > 0$. This allows us to sum over $k$ and conclude that
\[
\|S^2_m - S^{2,*}_m\|_{L^2 \to L^2}  \ \lesssim_P \ \sum_k 
\min({\mathcal A}_{p_{*}}(k), {\mathcal A}_{p_{*}}(k)^{-1})^{\epsilon_{*}} \ \lesssim_P \ 1,
\]
reducing matters to bounding $S^{2,*}_m$, uniformly in $m$ - in other words, we have safely removed term $\psi_{p_{*}}(y-t) s^{p_{*}}$ from the phase. 

We can now apply this argument iteratively, comparing $S^{2,*}_m$ to $S^{2,**}_m$ where
the phase in $S^{2,**}_m$ has both $\psi_{p_{*}}$ and $\psi_{p_{**}}$ removed and
$p_{*}, p_{**} \ge p_0$. Notice though that the same argument
above also allows us to remove an entire tail 
\[
{\tilde \psi}_{p_1}(y-t,s)  \ = \ \sum_{p\ge p_1} \psi_p(y-t) s^p \ \ {\rm for \ some} \ \ p_1 \ge p_0.
\]
In fact 
we may suppose that there is a $p_1 \ge p_0$ such that
$|\psi_{p_1}(y-t)| \sim_P |y-t|^{\ell_{1}}$ for some $\ell_{1} \ge 1$ and $\psi_p^{(\ell)}(0) = 0$
for all $p\ge p_1$ and all $\ell < \ell_{1}$. Otherwise $\psi_p \equiv 0$ for all $p \ge p_0$ and so 
${\tilde \psi}_{p_0} \equiv 0$.
Hence $|{\tilde \psi}_{p_1}(y-t,s)| \sim c_1 |(y-t)^{\ell_1} s^{p_1}|$ for some $c_1$ and so
${\tilde \psi}_{p_1}(y-t, s)$ can be treated in the same way as $\psi_{p_1}(y-t) s^{p_1}$ and thus be removed from the phase. The above iteration then removes the remaining terms with $p_0 \leq p_{*} < p_1$.

Hence for Theorem \ref{vw-analytic},  the uniform (in $m$) $L^2$ boundedness of $S^2_m$ is equivalent to the uniform 
(in $m$) $L^2$ boundedness of 
\[
H_m g(y) \ := \   \chi_m(y) \sum_{I\in {\mathcal F}^{0,m}} \int_{\mathbb R} \rho_I(\lambda, y, t)  \, g(t) \, dt
\]
where 
\[
\rho_I(\lambda, y, t)  \ =  \ \int_{\mathbb R} 
e^{2\pi i \lambda ( 2 y s^{p_0} + \phi(s) + \sum_{p=0}^{p_0 -1} \psi_p(y-t) s^p)} \,  \phi_I^{(I)}(s,y-t) \, ds.
\]
Note that $P_{p_0}(s,t) = \phi(s) + \sum_{0\le p < p_0} \psi_p(t) s^p$ is precisely the function
featuring in the statement of Theorem \ref{vw-analytic}.

\subsection{Back to the common analysis of Theorems \ref{main-uniform} and \ref{vw-analytic}}
To unify the notation somewhat, we will designate as ${\mathcal H}_m$ both the operator
\[
S^2_m g(y) \ = \ \chi_m(y) \sum_{I\in {\mathcal F}^{0,m}} \int_{\mathbb R} m_I(\lambda, y, t) \, g(t) \, dt
\]
when we refer to Theorem \ref{main-uniform} and the operator $H_m$ in the previous section
defined with $\rho_I$ instead of $m_I$ when we refer to Theorem \ref{vw-analytic}. Furthermore we relabel $\rho_I$
as $m_I$ so that when we refer to Theorem \ref{vw-analytic},
\[
m_I(\lambda, y, t)  \ =  \ \int_{\mathbb R} 
e^{2\pi i \lambda ( 2 y s^{p_0} + \phi(s) + \sum_{p=0}^{p_0 -1} \psi_p(y-t) s^p)} \,  \phi_I^{(I)}(s,y-t) \, ds
\]
and when we refer to Theorem \ref{main-uniform},
\[
m_I(\lambda, y, t)  \ =  \ \int_{\mathbb R} 
e^{2\pi i \lambda ( 2 y s + \phi(s) + \sum_{p\ge 0} \psi_p(y-t) s^p)} \,  \phi_I^{(I)}(s,y-t) \, ds.
\]
Of course the functions in the phase of $m_I$ are real-analytic for Theorem \ref{vw-analytic} and they
are polynomials for Theorem \ref{main-uniform}.

We split the operator ${\mathcal H}_m = {\mathcal H}_m^1  + {\mathcal H}_m^2 $ where
\[
{\mathcal H}_m^1 g(y) \ := \ 
 \chi_m(y)  \sum_{I \in {\mathcal F}^{0,m}_1}  \int_{\mathbb R}
m_I(\lambda, y, t) \,  g(t) \, dt, 
\]
with
\[
{\mathcal F}^{0,m}_1 = \{I = (j,k)\in {\mathcal F}^{0,m}: k \in K_1\} \  \ {\rm and}  \ \
K_1 = \{k : {\mathcal A}_p(k) \le 1, {\rm for \ all} \ p \in {\mathcal P}\}.
\]
The operator ${\mathcal H}_m^2$ is defined similarly
where the $k$ sum with $I \in {\mathcal F}^{0,m}_2$ is 
taken over the complementary set $K_2$ where at least one $p \in {\mathcal P}$ 
satisfies ${\mathcal A}_p(k) \ge 1$.

For ${\mathcal H}_m^1$, we proceed as in Section \ref{key-step}, using 
\eqref{Dk} to bound the difference ${\mathcal H}_m^1 - {\mathcal H}_m^{1,*}$ where ${\mathcal H}_m^{1,*}$ is defined
the same as ${\mathcal H}_m^1$ except with $m_I$ replaced by $m_I^{*}$ -- see \eqref{m*} (of course for
Theorem \ref{vw-analytic}, we need to adjust
appropriately the phase in $m_I^{*}$ -- we also note that the difference bound
\eqref{Dk} still holds for $m_I - m_I^{*}$ in the context of Theorem \ref{vw-analytic}). Hence
\eqref{Dk} implies that
\[
\|{\mathcal H}_m^1 - {\mathcal H}_m^{1,*}\|_{L^2 \to L^2} \ \lesssim \
\sum_{k\in K_1} {\mathcal A}_{p_{*}}(k)^{\epsilon_{*}} \ \lesssim
\ 1.
\]

Proceeding iteratively, we see that the uniform
boundedness of ${\mathcal H}_m^1$ is reduced to the uniform boundedness of
\[
L_m g(y) \ := \ \chi_m(y) \sum_{I  \in {\mathcal F}^{0,m}_1}  \int_{\mathbb R} \tau_I(\lambda, y, t) \,
e^{2\pi i \lambda \psi_0(y-t)} g(t) dt
\]
where
\[
\tau_I(\lambda, y, t) \ = \ \int_{\mathbb R} 
e^{2\pi i \lambda ( 2 y s^{p_0} + \phi(s))} \,  \phi_I^{(I)}(s,y-t) \, ds.
\]
We note that  $L_m g(y) = \chi_m(y) \int K_m(y, y-t) e^{2 \pi i \lambda \psi_0(y-t)} g(t) \, dt$ where
\[
K_m(y, y-t) \ = \ \sum_{k\in K_1} 2^{-k} K_m^{(k)}(y, 2^{-k}(y-t)) 
\]
and
\[
K_m^{(k)}(y,\tau) \ = \ \sum_{j: I = (j,k) \in {\mathcal F}^{0,m}_1} 
\int_{\mathbb R} e^{2\pi i \lambda (2y s^{p_0} + \phi(s))}  2^{-j}\phi_I (2^{-j}s, \tau) \, ds
\]
Hence $K_m$ is a variable Calder\'on-Zygmund kernel on ${\mathbb R}$; that is,
\begin{equation}\label{variable}
\int_{\mathbb R} K_m(y, \tau) \, d\tau \ = \ 0 \ \ {\rm for \ all} \ y  \ \ 
{\rm and} \ \ \forall \ell, \ |\partial^{\ell}_{\tau} K_m (y, \tau)| \lesssim |\tau|^{-\ell -1} \ 
\end{equation}
holds, uniformly in $m$ and $y$. This follows from an simple variant of
\eqref{sw-generalised}; more precisely, one sees that  \eqref{k} remains true with
$\phi_I$ replaced by any derivative $\partial^{(k)}_{t} \phi_I(s,t)$. 

This puts us in a position to appeal to a theorem of Ricci and Stein in \cite{RS1} on uniform $L^2$ bounds for
oscillatory singular integral operators
\[
T_{\lambda} g(y) \ = \ \int_{\mathbb R} K(y-t) e^{i\lambda \psi_o(y-t)} \, g(t) \, dt.
\]
When $\psi_0$ is a polynomial (which is the case for Theorem \ref{main-uniform}), 
Ricci and Stein establish $L^2$ bounds which are uniform in $\lambda$,
the Calder\'on-Zgymund kernel $K$ and the coefficients of $\psi_0$. In \cite{Pan}, Pan extended
this result to real-analytic phases $\psi_0$ (the case for Theorem \ref{vw-analytic}). 
Although their results are stated and proved for classical CZ kernels,
an examination of their arguments shows that the same results hold for variable CZ kernels described above in \eqref{variable}.
At the heart of their argument is a $T_{\lambda}^{*}T_{\lambda}$ argument applied to dyadic pieces
of the operator. Fortunately the order of the composition is immaterial (in fact they chose the order  
$T_{\lambda}^{*}T_{\lambda}$)
but for our variable CZ kernel $K_m$ above, it is important to take the order $T_{\lambda} T_{\lambda}^{*}$
so that the variable $y$ in the first argument of $K_m(y, y-t)$ does not interact with the integration defining
the kernels of the various $T_{\lambda} T_{\lambda}^{*}$s. 
We leave the details to the reader. This completes the analysis for the ${\mathcal H}_m^1$; they
define uniformly bounded $L^2$ operators.
  
For ${\mathcal H}^2_m$, our
goal will be to establish uniform $L^2$ bounds for the difference ${\mathcal H}_m^2 - T_m'$ where $T_m'$
is defined exactly the same as ${\mathcal H}_m^2$ except that $m_I(\lambda, y, t)$ is replaced by
\[
e_I(\lambda, y - t)  \ =  \ \int_{\mathbb R} 
e^{2\pi i \lambda ( \phi(s) + \sum_{p\ge 0} \psi_p(y-t) s^p)} \,  \phi_I^{(I)}(s,y-t) \, ds
\]
for Theorem \ref{main-uniform} and
\[
e_I(\lambda, y - t)  \ =  \ \int_{\mathbb R} 
e^{2\pi i \lambda ( \phi(s) + \sum_{0 \le p < p_0} \psi_p(y-t) s^p)} \,  \phi_I^{(I)}(s,y-t) \, ds
\]
for Theorem \ref{vw-analytic}. That is, for $\mathcal{H}^2_m$ we plan to remove the term $2ys^{p_0}$ from the phase this time. Note that the phase in the first integral is precisely
the original $P(s,y-t)$.

It is a simple matter to see that uniform boundedness of the family $\{T_m'\}$ is equivalent
to the uniform boundedness of the euclidean translation-invariant family $\{T_m\}$ where
\[
T_m g(y) \ = \ K_m * g(y) \ \ {\rm and} \ \ K_m(\tau) \ = \ \sum_{I\in {\mathcal F}^{0,m}_2} e_I(\lambda, \tau);
\]
thus $T_m$ is the same at $T'_m$ without the $\chi_m(y)$ factor in front.

In fact from the pointwise bound $|T_m' g(y)| \le |T_m g(y)|$, one direction is clear. Suppose
now that the family $\{T_m'\}$ is uniformly bounded in $L^2$ and decompose an
$L^2(\mathbb R)$ function $g = \sum_{\ell} g_{\ell}$ so that the support of 
${\tilde g}_{\ell}(t) := g_{\ell}(\ell 2^m + t)$ is contained in $\{|t| \sim 2^m\}$. 
Since for $I = (j,k) \in {\mathcal F}^{0,m}_2, \, k \ll m$, we see that if $|y-t| \sim 2^k$ and $|t| \sim 2^m$,
then $|y| \sim 2^m$ and so
\[
T_m g_{\ell} (y + \ell 2^m) \ = \ T_m {\tilde g}_{\ell} (y) \ = \ \chi_m(y) T_m {\tilde g}_{\ell} (y)
\ = \ T_m' {\tilde g}_{\ell} (y).
\]
Therefore, by almost disjointness of the supports,
\[
\|T_m g\|_{L^2}^2 \lesssim \sum_{\ell} \|T_m g_{\ell} \|_{L^2}^2 = \sum_{\ell} \|T_m' {\tilde g}_{\ell} \|_{L^2}^2
\lesssim \sum_{\ell} \|{\tilde g}_{\ell} \|_{L^2}^2 = \sum_{\ell} \|g_{\ell}\|_{L^2}^2 = \|g\|_{L^2}^2.
\]

The difference ${\mathcal H}_m^2 - T_m'$ is
\[
({\mathcal H}_m^2 - T_m' ) g(y) = \chi_m(y) \int_{\mathbb R} \sum_{I \in {\mathcal F}^{0,m}_2} [m_I(\lambda, y, t) -
e_I(\lambda, y-t)] \, g(t) \, dt
\]
and so we concentrate on bounding the difference
\[
{\mathcal D} \ := \ \sum_{ I \in {\mathcal F}^{0,m}_2} [m_I - e_I] \ = \ 
\sum_{k\in K_2} \sum_{j: I = (j,k) \in {\mathcal F}^{0,m}_2}
 [m_I - e_I] \ =: \ \sum_{k\in K_2} {\mathcal D}_k.
\]
We split $K_2  = \bigcup_{p\in {\mathcal P}} K_{2,p}$ where 
\[
K_{2,p} \ := \ \bigl\{ k \in K_2 :  {\mathcal A}_p(k) \ge {\mathcal A}_{p'}(k), \ \forall \ p' \in {\mathcal P}\bigr\}
\]
so that when $k \in K_{2,p}$, we have ${\mathcal A}_p(k) \ge 1$ (by definition of $K_1$).
This gives a corresponding
splitting of ${\mathcal H}_m^2 - T_m' = \sum_{p\in{\mathcal P}} ({\mathcal H}_m^2 - T_m')_p$ 
where the $I = (j,k) \in {\mathcal F}^{0,m}_2$
is restricted to $k \in K_{2,p}$. 

We claim that for $k\in K_{2,p}$,
\begin{equation}\label{Dkp}
|{\mathcal D}_k| \ \lesssim_P \ {\mathcal A}_p(k)^{-\epsilon_p} \ 2^{-k} \, \chi_{|y-t| \sim 2^k}
\end{equation}
for some $\epsilon_p>0$. If this is the case, then we have
\[
\|({\mathcal H}_m^2 - T_m')_p\|_{L^2 \to L^2} \ \lesssim_{p} \ 
\sum_{k\in K_{2,p}} {\mathcal A}_p(k)^{-\epsilon_p} \ \lesssim \ 1
\]
and so summing over $p \in {\mathcal P}$ gives the desired uniform bound for ${\mathcal H}_m^2 - T_m'$. 

To prove
\eqref{Dkp}, we fix $p$ and $k \in K_{2,p}$ and split
\[
{\mathcal D}_k \ = \sum_{j\in J_1} [m_I - e_I] \ + \ \sum_{j\in J_2} [m_I - e_I] \ := \ {\mathcal D}_k^1 +
{\mathcal D}_k^2
\]
into two parts; here $J_1 = \{ j : 2^{j } \le (|\lambda| 2^m)^{-1/p_0} {\mathcal A}_p (k)^{-\sigma_p}  \}$
for some $\sigma_p > 0$ and $J_2$ is the complementary range.

For ${\mathcal D}_k^1$, we use the difference in the phases of $m_I$ and $e_I$ to see  that
$$
|m_I(\lambda, y, t) - e_I(\lambda, y-t)| \ \lesssim \ |\lambda y| 2^{j p_0} \ \sim \ |\lambda| 2^m 2^{j p_0}
$$
and so
$$
|{\mathcal D}_k^1| \ \lesssim \ |\lambda| 2^m \, \Bigl[\sum_{j\in J_1} 2^{j p_0}\Bigr] \, 2^{-k} \chi_{|y-t| \sim 2^k}
\ \lesssim \ 
{\mathcal A}_p(k)^{-p_0 \sigma_p} \, \ 2^{-k} \chi_{|y-t| \sim 2^k},
$$
establishing \eqref{Dkp} for ${\mathcal D}_k^1$. For ${\mathcal D}_k^2$ we treat the terms
$m_I$ and $e_I$ separately, bounding 
$|{\mathcal D}_k^2| \le \sum_{j\in J_2} |m_I| + \sum_{j\in J_2} |e_I|$.

We will apply both \eqref{Ij-poly} and \eqref{Ij-analytic-p} to each $m_I$ and $e_I$ separately. The phase in $m_I$ is
$$
 \lambda \Bigl[ (2 y s^{p_0} + \sum_{p\ge 1} c_{p,0} s^p +  \sum_{p \ge 0}  \psi_p(y-t) s^p \Bigr]
$$
for Theorem \ref{main-uniform} whereas for Theorem \ref{vw-analytic}, the sum 
$\sum_{p =0}^{p_0 -1}  \psi_p(y-t) s^p$ is truncated. The
phase in $e_I$ is the same except the term $2 y s^{p_0}$ is not present. They both have the $s^p$ coefficient
$h_{p} := \lambda (c_{p,0} + \psi_{p}(y-t))$ unless $p=1$ and we are in the setting of Theorem \ref{main-uniform}.
Setting this case aside for the moment,
we apply \eqref{Ij-poly} and \eqref{Ij-analytic-p} to each $m_I$ and
$e_I$ with respect to this common coefficient $h_p$. Since for some $\ell_p \ge 1$,
 $|\psi_p(y-t)| \sim c_p 2^{\ell_p k}$ when $|y-t| \sim 2^k$, 
we see that there are only $O(1)$ values of $k$
where the bound $|h_{p}| \sim |\lambda| c_p 2^{\ell_p k}$ does not hold. We add these values
to the expectional set ${\mathcal E}$. 
Hence in this case, \eqref{Ij-poly} and \eqref{Ij-analytic-p} imply
\begin{equation}\label{m-e}
|m_I|, |e_I| \ \lesssim \ \bigl[|\lambda| c_p 2^{\ell_p k} 2^{p j}\bigr]^{-\epsilon_0} \ 2^{-k} \chi_{|y-t| \sim 2^k}
\end{equation}
for some $\epsilon_0 > 0$. 

If in the context of Theorem \ref{main-uniform} (so that $p_0 = 1$ and hence
the coefficient $c_{1,0}$ in $\phi(s)$ is zero) we are considering 
the case $p=1$, observe that the coefficient of $s$ for $m_I$, which is $h_1 = \lambda (2y + \psi_1(y -t))$, is different
from the coefficient of $s$ for $e_I$, $h_1 = \lambda \psi_1(y-t)$. However in both cases, except for a few
values of $k$ (which we toss into ${\mathcal E}$), we have $|h_1| \gtrsim |\lambda| c_1 2^{\ell_1 k}$
and so the estimate \eqref{m-e} holds in this case as well if one chooses $\sigma_1$ so that $0 < \sigma_1 < 1 $.

Summing the estimates \eqref{m-e}
over $j\in J_2$ establishes \eqref{Dkp} for ${\mathcal D}_k^2$ and hence ${\mathcal D}_k$. This shows
that the uniform $L^2$ boundedness of ${\mathcal H}_m$ is equivalent to the uniform
$L^2$ boundedness of $T_m$. 

Putting everything together, we see that the $L^2$ boundedness
of the original convolution operator $T_{P, {\mathcal F}}$ on the Heisenberg group ${\mathbb H}^1$ is equivalent
to the uniform in $m$ (and $\lambda$) $L^2$ boundedness of the euclidean convolution operators 
$T_m$. Recall the definition of the operators $T_m$ differ depending on whether
we are in the context of Theorem \ref{main-uniform} or Theorem \ref{vw-analytic}.
In the context of Theorem \ref{main-uniform}, the multiplier for $T_m$ is
\[
\int_{\mathbb R} K_m(t) e^{2\pi i \eta t} \, dt \ = \ \sum_{I\in {\mathcal F}^{0,m}_2}
\int\!\!\!\int_{{\mathbb R}^2} e^{2\pi i (\eta t + \lambda P(s,t))}
 \phi_I^{(I)} (s,t) \, ds dt
\]
and so the uniform $L^2$ boundedness of the $T_m$ is equivalent to showing that the above sum
of integrals is bounded
uniformly in the parameters $m$, $\lambda$ and $\eta$.

In the context of Theorem \ref{vw-analytic}, the
multiplier for $T_m$ is
\[
{\widehat{K_m^{\lambda}}}(\eta) \ = 
\ \int_{\mathbb R} K_m(t) e^{2\pi i \eta t} \, dt \ = \ \sum_{I\in {\mathcal F}^{0,m}_2}
\int\!\!\!\int_{{\mathbb R}^2} e^{2\pi i (\eta t + \lambda P_{p_0}(s,t))} \phi_I^{(I)} (s,t) \, ds dt
\]
and uniform boundedness is equivalent to showing that ${\widehat{K_m^{\lambda}}}(\eta)$ 
is uniformly bounded in $m$, $\lambda$ and $\eta$.

\section{The conclusion of the proof of Theorem \ref{vw-analytic}}
 
Consider the following truncations of the multiparameter singular Radon transform $R_{P_{p_0}, K}$ (from the
statement of Theorem \ref{vw-analytic}):
\[
R_{P_{p_0}, K_m} f(x,y) \ = \ \int\!\!\!\int_{{\mathbb R}^2} f(x - t, y - P_{p_0}(s,t)) \,  {K}_m (s,t) ds dt
\]
where 
\[
{K}_m(s,t) \ = \ \sum_{I \in {\mathcal F}^{0,m}_2} \phi_I^{(I)}(s,t)
\]
is a truncation of the product kernel $K$. The multiplier $M_m(\eta, \lambda)$ of
$R_{P_{p_0}, K_m}$ is precisely equal to
${\widehat{K_m^{\lambda}}}(\eta)$ above.

Thus the $L^2({\mathbb H}^1)$ boundedness of  $T_{P, {\mathcal F}}$ is equivalent to the uniform 
$L^2({\mathbb R}^2)$ boundedness of the truncations $R_{P_{p_0}, K_m}$ as stated in 
Theorem \ref{vw-analytic}. When $K(s,t) = {\mathcal K}(s,t) = 1/s t$ is the double Hilbert transform
kernel, the operator $R_{P_{p_0}, K}$ and its generalisations have been 
thoroughly investigated in several papers; see for example, \cite{CWW2}, \cite{CWW1}, \cite{P}
and \cite{PY}. In \cite{CWW2} it is shown
that $R_{P_{p_0}, {\mathcal K}}$ is bounded on $L^2$ if and only if every vertex of the Newton
diagram of $P_{p_0}$ has at least one even component. It is straightforward to check that the same
conclusion holds for the truncated operators $R_{P_{p_0}, {\mathcal K}_m}$

This completes the proof of Theorem \ref{vw-analytic}.

\section{The conclusion of the proof of Theorem \ref{main-uniform}}\label{main}

Consider the following truncations of the multiparameter singular Radon transform $R_{P, K}$ (from the
statement of Theorem \ref{main-uniform}):
\[
R_{P, K_m} f(x,y) \ = \ \int\!\!\!\int_{{\mathbb R}^2} f(x - t, y - P(s,t)) \,  {K}_m (s,t) ds dt
\]
where 
\[
{K}_m(s,t) \ = \ \sum_{I \in {\mathcal F}^{0,m}_2} \phi_I^{(I)}(s,t)
\]
is a truncation of the product kernel $K$. The multiplier $M_m(\eta, \lambda)$ of
$R_{P, K_m}$ is precisely equal to the multiplier of $T_m$; that is, 
\[
M_m(\lambda, \eta) \ = \ \sum_{I\in {\mathcal F}^{0,m}_2}
\int\!\!\!\int_{{\mathbb R}^2} e^{2\pi i (\eta t + \lambda P(s,t))}
 \phi_I^{(I)} (s,t) \, ds dt.
\]

Thus the uniform $L^2({\mathbb H}^1)$ boundedness of  $T_{P, {\mathcal F}}$ (where we seek
uniformity over $P \in {\mathcal V}_{\Delta}$ and the truncations ${\mathcal F}$)
is equivalent to the uniform 
$L^2({\mathbb R}^2)$ boundedness of $R_{P, K_m}$ where uniformity in $m$ is also
required. This is the main statement in Theorem \ref{main-uniform}. 
When $K(s,t) = {\mathcal K}(s,t) =  1/s t$ is the double Hilbert transform
kernel, we can apply Theorem 5.1 form \cite{RS3} exactly as we did for the Ricci-Stein theorem
from the Introduction to conclude that
$$
\sup_m \sup_{P\in {\mathcal V}_{\Delta}} \|R_{P, {\mathcal K}_m}\|_{L^2({\mathbb R}^2) \to L^2({\mathbb R}^2)}
\ < \ \infty
$$
if and only if every $\alpha = (\alpha_1, \alpha_2) \in \Delta$ has at least one even component. The 
{\it only if} part of the statement is an easy computation of the multiplier $M_m(\lambda, \eta)$ associated
to a single monomial $P(s,t) = s^j t^k$ where both $j$ and $k$ are odd (see \cite{Fef}). 

This completes the proof of Theorem \ref{main-uniform}.

\section{Appendix - proof of \eqref{Ij-analytic-p}}

In this appendix we give a proof of the oscillatory integral estimate \eqref{Ij-analytic-p}. Recall 
\[
I_j \ = \ \int_{\mathbb R} e^{2\pi i Q(s)} \, 2^{-j} \Phi(2^{-j} s) \, ds
\]
where $Q$ is either
\[
 \lambda \Bigl[ 2 y s^{p_0} + \sum_{p\ge 1} c_{p,0} s^p +  \sum_{1\le p < p_0}  \psi_p(y-t) s^p \Bigr]
\]
or
\[
 \lambda \Bigl[ \sum_{p\ge 1} c_{p,0} s^p +  \sum_{1\le p < p_0}  \psi_p(y-t) s^p \Bigr].
\]
For ease in notation, we will assume $Q$ is the latter. When considering the former instead, without loss of generality one may assume there
exists an $c_{p,0} \neq 0$ for some $p> p_0$; otherwise, we would be in the polynomial case where we can appeal to \eqref{Ij-poly}.

Let $p_n < \cdots < p_1$ enumerate the values of $1\le p < p_0$ such that $\psi_p \not\equiv 0$. In
this case, for each $1\le r \le n$, there is an $\ell_r \ge 1$ such that $|\psi_{p_r}(y-t)| \sim 2^{\ell_r k}$
whenever $|y-t| \sim 2^k$. Hence $h_r := \lambda ( c_{p_r, 0} + \psi_{p_r}(y-t))$ satisfies
$|h_r| \gtrsim |\lambda| 2^{\ell_r k}$ whenever $|y-t| \sim 2^k \ll 1$ and with this notation, \eqref{Ij-analytic-p}
reads 
\begin{equation}\label{Ij-p}
|I_j| \ \lesssim_{P} \ \bigl[|h_{p_r}| 2^{p j}\bigr]^{-\epsilon_r} \ \lesssim \ 
\bigl[|\lambda| 2^{\ell_r k} 2^{p_r j}\bigr]^{-\epsilon_r}
\end{equation}
for every $1\le r \le p_0 -1$ and for some $\epsilon_r> 0$.

We fix an $1\le L < p_0$ and establish \eqref{Ij-p} with $r = L$. First of all, we have $|\psi_{p_r}(y-t)s^{p_r}| \sim 2^{\ell_r k} 2^{p_r j}$ and thus let us name these quantities $\theta_r(k,j) : = 2^{\ell_r k} 2^{p_r j}$; they will be used to control the contribution of each term of $Q$ to some derivative of $Q$ itself.\\
We introduce a sequence of small parameters $0 < \delta_1 \ll \delta_2 \ll \cdots \ll \delta_{L-1} \ll 1$ depending on $P$, which will be chosen later, and define for each $1\leq r \leq L$ sets 
\begin{align*}
R_r := \Big\{ j \; : \; \theta_1(k,j) &< \delta_1 \theta_L(k,j), \\
& \vdots \\
\theta_{r-1}(k,j) &< \delta_{r-1} \theta_L(k,j), \\
& \text{ and } \\
\theta_r(k,j) &\geq  \delta_r \theta_L(k,j)\Big\}.
\end{align*}
Notice that for $R_1$ the first conditions are vacuous and we only stipulate $\theta_1(k,j) \geq \delta_1\theta_L(k,j)$, and for $R_L$ the last condition is vacuous and we only stipulate $\theta_s(k,j) < \delta_s \theta_L(k,j)$ for all $s =1,\ldots, L-1$. It is immediate to see that these sets form a partition of the set of all possible $j$'s.\\
Suppose that $j \in R_r$ for some $1\le r \le L$. We examine the $p_r$-th derivative of $Q$:
\[
Q^{(p_r)}(s) = \lambda \Big[\sum_{i=1}^r \frac{p_i!}{(p_i - p_r)!} \psi_{p_i}(y-t) s^{p_i - p_r} +
 \frac{p_{*}!}{(p_{*} - p_r)!} c_{*} s^{p_{*} - p_r} + O(s^{p_{*} - p_r +1}) \Big]
\]
 where $p_{*} \ge p_r$ is the first exponent such that $c_{p_{*},0} \neq 0$. Noting $|s| \sim 2^j \ll 1$ and $j \in R_r$, the contribution of the mixed terms with $i<r$ is at most
\[ C_P \sum_{i<r} \theta_i(k,j)2^{-p_r j} \leq C_P(\delta_1 + \ldots + \delta_{r-1}) \theta_L(k,j)2^{-p_r j}, \]
while the contribution of the mixed term with $i=r$ is $\sim \theta_r(k,j)2^{-p_r j} > \delta_r \theta_L(k,j)2^{-p_r j}$. By choosing the constants $\delta_i$ to be sufficiently small (depending on $P$) and decreasing fast enough we have then
\[ \Big|\sum_{i=1}^r \frac{p_i!}{(p_i - p_r)!} \psi_{p_i}(y-t) s^{p_i - p_r}\Big| \gtrsim  \theta_L(k,j)2^{-p_r j}. \]
As for the contribution of the remaining terms, we have $|c_{*}| s^{p_{*} - p_r}\sim 2^{(p_{*} - p_r)j}$. Hence if $\theta_L(k,j) \not\sim 2^{p_{*}j}$ we have $|Q^{(p_r)}(s)|\gtrsim |\lambda|\theta_L(k,j)2^{-p_r j}$, implying that
\[
 |I_j| \ \lesssim \ \bigl( 2^{p_r j} |\lambda| \theta_L(k,j)2^{-p_r j} \bigr)^{-1/p_r} = \bigl(|\lambda| \theta_L(k,j)\bigr)^{-1/p_r}
\]
by van der Corput's lemma. Hence \eqref{Ij-p} holds in this case.

As for the case $2^{p_{*}j} \sim \theta_L(k,j)$, we have the bound
$|Q^{(p_{*})}(s)| \gtrsim 1$ since every $2^{\ell_r k} \ll 1$. Another application of van der Corput's lemma shows
\[
|I_j| \lesssim \bigl( 2^{p_{*} j} |\lambda| \bigr)^{-1/p_{*}}  \sim \bigl(|\lambda| \theta_L(k,j)\bigr)^{-1/p_{*}}.
\]
This completes the proof of \eqref{Ij-p}.

\bibliography{multipar_sing_int_H_bibliography}
\bibliographystyle{amsplain}
\end{document}